\documentclass{amsart}
\usepackage{amscd}
\usepackage{amsmath}
\usepackage{amssymb}
\newtheorem{theorem}{Theorem}[section]
\newtheorem{lemma}[theorem]{Lemma}
\newtheorem{cor}[theorem]{Corollary}
\newtheorem{prop}[theorem]{Proposition}

\theoremstyle{remark}

\numberwithin{equation}{section}

\begin{document}
\title{Subshifts and $C^*$-algebras from one-counter codes}
\author{Wolfgang Krieger}
\address{Institute for Applied  Mathematics,
University of Heidelberg,
Im Neuenheimer Feld 294,
69120 Heidelberg, Germany
}
\email{krieger@math.uni-heidelberg.de}
\thanks{
This work  was supported  by JSPS Grant-in-Aid for Scientific Reserch 
(N0.\ 20540215).}

\author{Kengo Matsumoto}
\address{Department of Mathematical Sciences, 
Yokohama City  University,
22-2 Seto, Kanazawa-ku, Yokohama, 236-0027 Japan}
\email{kengo@yokohama-cu.ac.jp}



\begin{abstract}
We introduce a class of subshifts  under the name of \lq\lq standard  
one-counter shifts \rq\rq. The standard one-counter shifts are the Markov  
coded systems of certain Markov codes that belong to the family of  
one-counter languages. We study topological conjugacy and flow  
equivalence of standard one-counter shifts. To subshifts there are  
associated C*-algebras by their $\lambda$-graph systems. We describe a  
class of standard one-counter shifts with the property that the  
C*-algebra associated to them is simple, while  the C*-algebra that is  
associated to their inverse is not. This gives examples of subshifts  
that  are not flow equivalent to their inverse. For a family of highly  
structured standard one-counter shifts we compute the K-groups.

\end{abstract}
\maketitle
\def\Zp{{ {\Bbb Z}_+ }}
\def\C{{{\mathcal C}}}
\def\V{{{\mathcal V}}}
\def\PL{{ {}^X\!{\frak L} }}
\def\FL{{ {\frak L}^X }}
\def\det{{{\operatorname{det}}}}
\def\trace{{{\operatorname{trace}}}}
\def\card{{{\operatorname{card}}}}
\def\M{{ {\mathcal M} }}
\def\H{{ {\mathcal H} }}
\def\K{{ {\mathcal K} }}
\def\CNC{{ {\mathcal C}^{(N)}_{counter} }}
\def\CNR{{ {\mathcal C}^{(N)}_{reset} }}
\def\sccnc{{ sc({\mathcal C}^{(N)}_{counter}) }}
\def\sccnr{{ sc({\mathcal C}^{(N)}_{reset}) }}
\def\sccnrt{{ sc(({\mathcal C}^{(N)}_{reset})^{rev}) }}
\def\Hom{{{\operatorname{Hom}}}}
\def\Ker{{{\operatorname{Ker}}}}
\def\id{{{\operatorname{id}}}}
\def\Ext{{{\operatorname{Ext}}}}
\def\OCNC{{ {\mathcal O}_{sc({\mathcal C}^{(N)}_{counter})   } }}
\def\OCNR{{ {\mathcal O}_{sc({\mathcal C}^{(N)}_{reset}) } }}
\def\OCRT{{ {\mathcal O}_{sc(({\mathcal C}^{(N)}_{reset})^{rev}) } }}
\def\synsyn{{ \Sigma_{synchro}(X) }}
\def\OPR{{ \Omega^+_{reset}(X) }}          
\def\OMR{{ \Omega^-_{reset}(X) }}
\def\OPC{{ \Omega^+_{counter}(X) }}          
\def\OMC{{ \Omega^-_{counter}(X) }}

\bigskip

\section{Introduction}
Let $\Sigma$ be a finite alphabet.
We use notation like
$$
x_{[i,k]} = {(x_j)}_{i\le j \le k}, \qquad
x \in \Sigma^{\Bbb Z},\, \,  i,k \in {\Bbb Z},
\, \, i \le k,  
$$
and we denote by 
$
x_{[i,k]}
$
also the word that is carried by the block
$
x_{[i,k]}.
$
The length of a word $a$ is denoted by $\ell(a)$.
On the shift space ${\Sigma}^{\Bbb Z}$ there acts the shift by
$$ 
x \longrightarrow
{(x_{i+1})}_{i \in {\Bbb Z}},
\qquad
x = {(x_i)}_{i \in {\Bbb Z}} \in \Sigma^{\Bbb Z}.
$$
A closed shift-invariant subset of 
$\Sigma^{\Bbb Z}$ is called a subshift.
For an introduction to the theory of 
subshifts  see \cite{Ki, LM}.
A word is called admissible for a subshift
if it appears in a point of the subshift.
We denote the language of admissible words of a subshift 
$X\subset \Sigma^{\Bbb Z}$
by ${\mathcal L}(X)$ 
and set  
${\mathcal L}_n(X) = \{ a \in {\mathcal L}(X)\mid \ell(a) = n \}, n \in {\Bbb N}.
$
A subshift
$X\subset \Sigma^{\Bbb Z}$
is uniquely determined by 
${\mathcal L}(X)$.
For a subshift
$X\subset \Sigma^{\Bbb Z}$
and for $I_{-}, I_+ \in {\Bbb Z}, I_{-} < I_{+}$,
one has a topological conjugacy
$$
x \longrightarrow (x_{[i+I_{-}, i+I_{+}]})_{i \in {\Bbb Z}}, \qquad (x \in X)
$$
of $X$ onto the higher block system 
$X^{\langle[I_{-},I_{+}]\rangle}$ of $X$.

Among the first  
 examples of subshifts are 
the topological Markov shifts.
Using a matrix 
$(A(\sigma,\sigma'))_{\sigma, \sigma' \in \Sigma},$
$$
A(\sigma,\sigma') \in \{0,1 \},\qquad  \sigma, \sigma' \in \Sigma,
$$
that has in every row and every column at least one entry that is equal to $1$
as a transition matrix one obtains a topological Markov shift
$tM(\Sigma,A)$ 
by setting
$$ 
tM(\Sigma,A) = 
 \{ {(\sigma_i)}_{i \in {\Bbb Z}} \in \Sigma^{\Bbb Z}
\mid A(\sigma_i, \sigma_{i+1}) =1, i \in \Bbb Z \}.
$$
For $n >1$ 
the $n$-block system 
$(\Sigma^{\Bbb Z})^{\langle[1,n]\rangle}$
of the shift 
on $\Sigma^{\Bbb Z}$ is a topological Markov shift 
with a transition
matrix
$A^{(n)}$ that is given by
$$
A^{(n)}(a, a')=
\begin{cases}
1 & \text{if } a_{(1,n]} = a'_{[1,n)},\\
0 & \text{if } a_{(1,n]} \ne a'_{[1,n)}
\end{cases},
\qquad a,a' \in {\Sigma}^n. 
$$
A subshift
$ X \subset {\Sigma}^{\Bbb Z}$ is said to be 
of finite type if there is a finite set 
${\frak F}$ of words 
in the alphabet $\Sigma$
such that
$(\sigma_i)_{i \in {\Bbb Z}} \in X$
precisely if
no  word in
${\frak F}$
appears in 
$(\sigma_i)_{i \in {\Bbb Z}}$.
A subshift is topologically conjugate to a subshift 
of finite type if and only if it is of finite type
\cite{Ki, LM}.
We formulate this theorem equivalently as:
\begin{theorem}
Let $X \subset \Sigma^{\Bbb Z}$ be a subshift that is 
topologically conjugate to a topological Markov shift.
Then there exists an  $n_\circ \in {\Bbb N}$
such that 
\begin{equation*}
X^{{\langle [1,n]\rangle}} = tM({\mathcal L}_n(X),(A^{(n)}(a,a'))_{a,a' \in {\mathcal L}_n(X)}), \qquad n \ge n_\circ.
\end{equation*}
\end{theorem}

The coded system \cite{BH} of a formal language ${\mathcal C}$
in a finite alphabet
$\Sigma$
is the subshift
that is obtained as the closure of the set of points
in $\Sigma^{\Bbb Z}$
that carry bi-infinite concatenations of words in
${\mathcal C}$.
${\mathcal C}$ can here always be chosen 
to be a prefix code.
The property of being coded 
is an invariant of topological conjugacy.
We denote the coded system
of a code 
${\mathcal C}$  by
$sc({\mathcal C})$.
More generally a Morkov code (see \cite{Ke}) 
is given by a formal language  ${\mathcal C}$
of words in a finite alphabet $\Sigma$ together with
a  finite index set $\Gamma$
and, mappings 
$s: {\mathcal C} \longrightarrow \Gamma,
 t: {\mathcal C} \longrightarrow \Gamma
$
and a transition matrix
$(A(\gamma, \gamma'))_{\gamma,\gamma' \in \Gamma}, A(\gamma, \gamma')\in \{0,1\}, 
\gamma,\gamma' \in \Gamma.
$
From a Markov code $({\mathcal C},s,t)$
one obtains the Markov coded system
$scM({\mathcal C})$ as the subshift that is the closure
of the set of points $x \in \Sigma^{\Bbb Z}$
such that there are indices 
$i_k \in {\Bbb Z}, k \in {\Bbb Z}$,
$i_k < i_{k+1}, k \in {\Bbb Z}$
such that 
$x_{[i_k,i_{k+1})} \in {\mathcal C}, k \in {\Bbb Z}$,
and such that 
$$
A(t(x_{[i_{k-1},i_{k})}), s(x_{[i_k,i_{k+1})}) ) = 1, \qquad k \in {\Bbb Z}.
$$
With the alphabet
$\{ a_n \mid 1 \le n \le N \} \cup\{\alpha_-, \alpha_+\}, N \in {\Bbb N}$,
consider the codes  
\begin{align*}
\CNR
 = & \{ \alpha_{-}^k \alpha_{+}^m a_n \mid 1 \le n \le N, m, k \in {\Bbb N}, m \le k \}\\
\intertext{ and with the alphabet $\{ b_n \mid 1 \le n \le N \} \cup\{\alpha_-, \alpha_+\}, N \in {\Bbb N}$,
consider the codes}
\CNC
& = \{ \alpha_{-}^k \alpha_{+}^k b_n \mid 1 \le n \le N, k \in {\Bbb N} \}.
\end{align*}
The  coded systems 
$\sccnr$, $ \sccnc $ 
and
$\sccnr\cup\sccnc$
serve us as prototypes for a class of subshifts that we will call
standard one-counter shifts.
(Compare here \cite[Example 1, p. 561]{Bl},  
\cite[Example II, p. 449]{Ke}, \cite[Example 6.1, p. 896]{KM2002}). 
We arrive at a description of this class of subshifts 
by observing the behavior of 
$\sccnr$, $\sccnc $ 
and of
$\sccnr\cup\sccnc$
and by abstracting the essential structural properties
that these coded systems are to share with the standard one-counter shifts. 
$\sccnr$ and 
$\sccnr\cup\sccnc$
are prototypes 
of what we will call standard one-counter shifts with reset.
To every standard one-counter shift 
$X \subset {\Sigma}^{\Bbb Z}$
there is associated  a unique Markov code
${\mathcal C}^{(X)}$ such that 
$X = scM({\mathcal C}^{(X)})$ 
and such that a version of Theorem 1.1 holds.
A formal language is called a one-counter language 
if it is recognized by a push down automaton 
with one stack symbol \cite{F,FMR,HU}.
The Markov code ${\mathcal C}^{(X)}$ 
that is associated to a standard one-counter shift 
$X \subset {\Sigma}^{\Bbb X}$
is a one-counter language.

Given a subshift $X \subset {\Sigma}^{\Bbb X}$
a word $v \in {\mathcal L}(X)$
is called synchronizing if for 
$u, w \in {\mathcal L}(X)$ 
such that 
$ uv, vw \in {\mathcal L}(X)$
also
$ u v w \in {\mathcal L}(X)$.
A topologically  transitive
subshift is called synchronizing 
if it has a synchronizing word.
Before turning in Section 3
to the standard one-counter codes
we introduce in Section 3 auxiliary notions
for synchronizing subshifts.
We introduce strongly synchronizing subshifts 
as the subshifts in which synchronizing symbols appear 
uniformly close to synchronizing words,
and we introduce sufficiently synchronizing subshifts 
as the subshifts that have a strongly  synchronizing higher block system.

$\lambda$-graph systems (as introduced in \cite{Ma1999c})
are labeled directed graphs that are equipped with a shift like map $\iota$.
A $\lambda$-graph system ${\frak L}$ gives rise to a $C^*$-algebra
${\mathcal O}_{\frak L}$.
To a subshift $X$ there is invariantly associated  
a future  $\lambda$-graph system 
${}^X\!{\frak L}$
that is based on the future equivalences of the pasts in
$X_{(-\infty,0]}$ (as in \cite{KM2002})
and 
there is invariantly associated  
a past  $\lambda$-graph system 
${\frak L}^X$
that is based on the past equivalences of the futures in
$X_{[0,-\infty)}$ (as in \cite{Ma1999c}).
The future 
and the past $\lambda$-graph systems of a subshift are
time symmetric to 
one-another:
the future $\lambda$-graph system
of a subshift is identical to 
the past $\lambda$-graph system of its inverse and vice versa.
For a standard one-counter shift $X$ we will see that 
${\mathcal O}_{{}^X\!{\frak L}}$
is simple if and only if $X$ has reset and that 
${\mathcal O}_{{\frak L}^X}$ is not simple.
Since the stable isomorphism class of
${\mathcal O}_{{\frak L}^X}$ is an invariant of flow equivalence \cite{Ma2001b},
a standard one-counter shift with reset is not flow equivalent to its inverse.
For the one-counter shifts 
$\sccnr, \sccnrt$, we will show that 
\begin{align*}
K_0( \OCRT ) & \cong K_0( \OCNR ) \cong {\Bbb Z}/N {\Bbb Z} \oplus {\Bbb Z},\\
K_1( \OCRT ) & \cong K_1( \OCNR ) \cong  0.
\end{align*} 
The one-counter code $\CNC$ is equal to its reversal $(\CNC)^{rev}$.
The K-groups of the $C^*$-algebra have been computed in \cite{KM4}
as
\begin{align*}
K_0(\OCNC ) & \cong {\Bbb Z}/N {\Bbb Z} \oplus {\Bbb Z}^2,\\
K_1(\OCNC ) & \cong {\Bbb Z}.
\end{align*}
For another computation of K-groups of one-counter shifts see \cite{Ma2,Ma16}.
\medskip

For a subshift $X \subset \Sigma^{\Bbb Z}$
we set 
\begin{equation*}
X_{[i,k]} = \{ x_{[i,k]} \mid x \in X\}, 
\quad \qquad i,k \in {\Bbb Z}, 
\quad i \le k.
\end{equation*}
We set also
\begin{align*}
\Gamma_k^+(a)
& = \{b \in X_{(n,n+k]} \mid (a,b) \in X_{[m,n+k]} \},\quad  k \in {\Bbb N},\\
\Gamma_\infty^+(a)
& = \{x^+ \in X_{(n,\infty)} \mid (a,x^+) \in X_{[m,\infty)} \}, \quad n, m  \in {\Bbb Z}, \, m < n, \, 
a \in X_{[m,n]}.
\end{align*}
$\Gamma^-$ has the time symmetric meaning.

We recall that, given subshifts $X \subset \Sigma^{\Bbb Z}$,
 $\tilde{X} \subset \tilde{\Sigma}^{\Bbb Z}$,
 and a topological conjugacy 
 $\tilde{\varphi}:\tilde{X} \longrightarrow X$, 
 there is for some $L \in \Zp$ a block mapping
 $\tilde{\varPhi}:\tilde{X}_{[-L,L]} \longrightarrow \Sigma$,
 such that 
 $$
\tilde{\varphi}(\tilde{x}) 
= (\tilde{\varPhi}(\tilde{x}_{[i-L,i+L]}))_{i \in {\Bbb Z}}, 
\qquad \tilde{x} \in \tilde{X}.
$$
We set 
$$
\tilde{\varPhi}(\tilde{a}) =
(\tilde{\varPhi}(\tilde{a}_{[j-L,j+L]}))_{i+L \le j \le k-L},
\qquad \tilde{a} \in \tilde{X}_{[i,k]}, \quad i, k \in {\Bbb Z}, 
\quad k-i >2L.
$$ 
We use similar notation for words.

\section{Strong synchronization}
The first lemma is well known.
We include the proof for completeness.

\begin{lemma}
Let 
$
\tilde{X} \subset \tilde{\Sigma}^{\Bbb Z},
X \subset \Sigma^{\Bbb Z}
$
be subshifts and let
$\varphi:\tilde{X} \longrightarrow X$
be a topological conjugacy.
Let
$\tilde{L}, L \in \Zp$
be such that 
$[-\tilde{L}, \tilde{L}]$ is a coding window for $\varphi$ and 
$[-L, L]$ is a coding window for $\varphi^{-1}$.
Let
$\tilde{x} \in \tilde{X}, x = \varphi(\tilde{x})
$
and
$I_{-}, I_{+} \in {\Bbb Z}, I_{-} \le I_{+}.$
Let 
$x_{[I_{-}, I_{+}]}$ be synchronizing.
Then
$\tilde{x}_{[I_{-}-\tilde{L}-L, I_{+}+ \tilde{L} + L]}$ is synchronizing.
\end{lemma}
\begin{proof}
Let
\begin{align*}
\tilde{y}^{-} & \in \Gamma^{-}_{\infty}(\tilde{x}_{[I_{-}-\tilde{L}-L, I_{+}+ \tilde{L} + L]}), \\
\tilde{y}^{+} & \in \Gamma^{+}_{\infty}(\tilde{x}_{[I_{-}-\tilde{L}-L, I_{+}+ \tilde{L} + L]}), 
\end{align*}
and let
$y^{-} \in \Gamma^{-}_{\infty}(x_{[I_{-}, I_{+}]}), 
 y^{+} \in \Gamma^{+}_{\infty}(x_{[I_{-}, I_{+}]}),
$
be given by 
\begin{align*}
\widetilde{\varPhi}(\tilde{y}^{-}, \tilde{x}_{[I_{-}-\tilde{L}-L, I_{+}+ \tilde{L} + L]}) & = (y^{-}, x_{[I_{-}, I_{+}]}),\\
\widetilde{\varPhi}( \tilde{x}_{[I_{-}-\tilde{L}-L, I_{+}+ \tilde{L} + L]},\tilde{y}^{+}) & = ( x_{[I_{-}, I_{+}]},y^{+}).
\end{align*}
One has 
$(y^{-}, x_{[I_{-}, I_{+}]}, y^{+}) \in X$ 
and
$$
\varphi^{-1}(y^{-}, x_{[I_{-}, I_{+}]}, y^{+}) 
=(\tilde{y}^{-}, \tilde{x}_{[I_{-}-\tilde{L}-L, I_{+}+ \tilde{L} + L]},\tilde{y}^{+})
$$
and the lemma follows.
\end{proof}
For a subshift $X \subset {\Sigma}^{\Bbb Z}$,
we denote 
the set of its synchronizing symbols by $\synsyn$.
\begin{lemma}
Let 
$
\tilde{X} \subset \tilde{\Sigma}^{\Bbb Z},
X \subset \Sigma^{\Bbb Z}
$
be subshifts and let a topological conjugacy
$\varphi:\tilde{X} \longrightarrow X$
be given by a one-block map
$\tilde{\varPhi}:\tilde{\Sigma}\longrightarrow \Sigma.$
Let
$ L \in \Zp$
be such that 
 $\varphi^{-1}$ has coding window
$[-L, L]$ and set
$\hat{\varPhi}(\tilde{x}_{[-L,L]}) =\tilde{\varPhi}(\tilde{x}_0),  
\tilde{x}_{[-L,L]} \in \tilde{X}_{[-L, L]}.$
Then 
$$
\hat{\varPhi}^{-1}(\synsyn) \subset \Sigma_{synchro}(\tilde{X}^{\langle[-L,L]\rangle}).
$$
\end{lemma}
\begin{proof}
Apply Lemma 2.1.
\end{proof}

We say that a synchronizing subshift $X \subset \Sigma^{\Bbb Z}$ is strongly synchronizing
if there exists a $Q \in \Zp$ such that the following holds:
if $ x \in X$ and $I_{-}, I_{+} \in {\Bbb Z}, I_{-} < I_{+}$
are such that  
$x_{[I_{-}, I_{+}]}$ is synchronizing, then 
there exists an index $i,$  $I_{-} - Q \le i \le I_{+} + Q$
such that 
$x_{i}$ is a synchronizing symbol.

The higher block systems of a strongly synchronizing subshift are also 
strongly synchronizing.
We say that a subshift $X \subset \Sigma^{\Bbb Z}$ is sufficiently  synchronizing
if it has strongly synchronizing higher block systems.

\begin{prop}
Sufficient synchronization is an invariant of the topological conjugacy of subshifts.
\end{prop}
\begin{proof}
To prove the lemma
it is by Lemma 2.2 enough to consider the case of subshifts 
$
X \subset \Sigma^{\Bbb Z},
\tilde{X} \subset \tilde{\Sigma}^{\Bbb Z}
$
and of a topological conjugacy 
$\varphi:X \longrightarrow \tilde{X}$
that is given by a one-block map 
$\varPhi:\Sigma \longrightarrow \tilde{\Sigma}$
such that
\begin{equation}
\varPhi^{-1}(\tilde{\Sigma}_{synchro}(\tilde{X}) \subset \synsyn
\label{eqn:1}
\end{equation}
with $\tilde{X}$ strongly synchronizing and to show that 
$X$ is strongly synchronizing.
Let $L \in \Zp$ be such that 
$\varphi^{-1}$ has the coding window
$[-L,L]$ and let
$\tilde{Q}\in \Zp$ be such that
for 
$
\tilde{x} \in \tilde{X},
\tilde{I}_{-}, \tilde{I}_{+} \in {\Bbb Z}, 
\tilde{I}_{-} < \tilde{I}_{+}
$
such that 
$\tilde{x}_{[\tilde{I}_{-}, \tilde{I}_{+}]}$
is synchronizing,
one has
an $\tilde{i},$  $\tilde{I}_{-} - Q \le \tilde{i} \le \tilde{I}_{+} + Q$
such that 
$\tilde{x}_{\tilde{i}}$ is a synchronizing symbol.
 Then one has
 for $ x \in X$, $I_{-}, I_{+} \in {\Bbb Z}, I_{-} < I_{+}$
such that  
$x_{[I_{-}, I_{+}]}$ is synchronizing, 
by Lemma 2.1 that
$\varphi(x)_{[I_{-},I_{+}]}$
is synchronizing.
It follows that  
there exists an  $i\in {\Bbb Z},$  
$I_{-} -L -\tilde{Q} \le i \le I_{+} + L + \tilde{Q}$
such that 
$\varPhi(x_{i})$ is a synchronizing.
By (\ref{eqn:1}),
then $x_{i}$ is synchronizing.
\end{proof}

\section{Standard one-counter shifts}
{\bf 3 a. The structure of standard one-counter shifts}

Let $X \subset \Sigma^{\Bbb Z}$ be a topologically transitive subshift.
We call a pair  
$((\alpha_{-})_{i\in {\Bbb Z}},
(\alpha_{+})_{i\in {\Bbb Z}})
$
of fixed points of $X$ a characteristic pair, 
if it is the unique pair of fixed points
 that satisfies the following conditions $(a), (b) $ and $(c^-)$,
and a condition $(c^+)$
that is symmetric to condition $(c^-)$:

$(a)$ $X$ has a unique  orbit $O_X$  
that contains all points that are 
left asymptotic to $(\alpha_{-})_{i\in {\Bbb Z}}$
and 
right asymptotic to $(\alpha_{+})_{i\in {\Bbb Z}}$,
and that do not contain a synchronizing word.

$(b)$ $X$ has a point that is
left asymptotic to $(\alpha_{+})_{i\in {\Bbb Z}}$
and 
right asymptotic to $(\alpha_{-})_{i\in {\Bbb Z}}$
and that contains a synchronizing word.

$(c^-)$  There exists a $K \in {\Bbb N}$ 
such that the following holds:
If $ x \in X$ and $I_{-}, I_{+} \in {\Bbb Z}, I_{-} \le I_{+},$
are such that  
$x$ is
right asymptotic to $(\alpha_-)_{i\in {\Bbb Z}}$,
and 
$x_{[I_{-}, I_{+}]}$ is synchronizing, 
and
$x_{(I_+, I_+ +k]}$,
is not synchronizing,
$ k \in {\Bbb N}$,
then 
there exists an index $i,$  $I_{-} < i \le I_{+} + K$,
such that 
$x_{j} = \alpha_{-}, j \ge i$.
\begin{prop}
Let $X \subset \Sigma^{\Bbb Z}$ be a topologically transitive subshift
with a characteristic pair 
$((\alpha_{-})_{i\in {\Bbb Z}},(\alpha_{+})_{i\in {\Bbb Z}})
$
of fixed points, 
and let
$\tilde{\varphi}$ be a topological conjugacy of a subshift 
$\tilde{X} \subset {\tilde \Sigma}^{\Bbb Z}$
onto $X$.
Then 
$({\tilde{\varphi}}^{-1}((\alpha_{-})_{i\in {\Bbb Z}}),
{\tilde{\varphi}}^{-1}((\alpha_{+})_{i\in {\Bbb Z}}))
$
is a characteristic pair of fixed points of $\tilde{X}$.
\end{prop}
\begin{proof}
Conditions $(a), (b)$, $(c^-), (c^+)$,
being satisfied by
$((\alpha_{-})_{i\in {\Bbb Z}},(\alpha_{+})_{i\in {\Bbb Z}}),
$
the proposition follows by means of Lemma 2.1. 
\end{proof}

We introduce notation that we use for a synchronizing subshift 
$X \subset {\Sigma}^{\Bbb Z}$, that has a characteristic pair 
$((\alpha_{-})_{i\in {\Bbb Z}},(\alpha_{+})_{i\in {\Bbb Z}})
$
of fixed points.
For  $\sigma_{-} \in \synsyn$
we denote by ${\mathcal D}(\sigma_{-},\alpha_{-})$ 
the set of words $d^- \in {\mathcal L}(X)$,
that do not contain a synchronizing symbol, and that do not end with $\alpha$,  
such that  
\begin{equation*}
\sigma_{-} d^- \in \bigcap_{k \in {\Bbb N}} \Gamma^-(\alpha_-^k), 
\end{equation*}
and for $\sigma_{-}^+ \in \synsyn$
we denote by ${\mathcal D}(\sigma_{-}^+,\alpha_{+})$ 
the set of words $d_+^- \in {\mathcal L}(X)$,
that do not contain a synchronizing symbol, and that do not begin with $\alpha$,  
such that 
\begin{equation*}
\sigma_{-}^+ d_+^- \in \bigcap_{k \in {\Bbb N}} \Gamma^-(\alpha_+^k). 
\end{equation*}
We set 
\begin{align*}
\Sigma_{-}(X) & = 
\{ \sigma_{-} \in \synsyn \mid {\mathcal D}(\sigma_{-},\alpha_{-}) \ne \emptyset \},\\
\Sigma_{-}^+(X) & = 
\{ \sigma_{-}^+ \in \synsyn \mid {\mathcal D}(\sigma_{-}^+,\alpha_{+}) \ne \emptyset \}.\\
\end{align*}
${\mathcal D}(\alpha_{+},\sigma_+,), {\mathcal D}(\alpha_-,\sigma_{+}^-,)$
and
$\Sigma_{+}(X),\Sigma_{+}^-(X)$
have the symmetric meaning.
\begin{lemma}
For a strongly synchronizing subshift
 $X \subset \Sigma^{\Bbb Z}$ that has a characteristic pair 
$((\alpha_{-})_{i\in {\Bbb Z}}, 
(\alpha_{+})_{i\in {\Bbb Z}})
$
of fixed points,
the sets $\Sigma_{-}(X)$ and $\Sigma_{+}(X)$
are not empty and the sets 
${\mathcal D}(\sigma_{-},\alpha_{-}), \, \sigma_{-}\in \Sigma_{-}(X)$
and
${\mathcal D}(\alpha_{+},\sigma_{+}), \, \sigma_{+}\in \Sigma_{+}(X)$
are finite.
\end{lemma}
\begin{proof}
We show that  $\Sigma_{-}(X)$ 
is not empty, and that the sets 
${\mathcal D}(\sigma_{-},\alpha_{-}), \, \sigma_{-}\in \Sigma_{-}(X)$
are finite.
By condition $(b)$ there exists
 an $x \in X$ that contains a synchronizing word and 
 that is right asymptotic to 
$(\alpha_{-})_{i \in {\Bbb Z}}$.
The assumption that $X$ is strongly synchronizing 
implies that $x$ contains a synchronizing symbol.
Let $i \in {\Bbb Z}$ be such that 
$x_i \in \synsyn, x_{i+K} \not\in \synsyn, K \in {\Bbb N}$.
If here $x_{i+K} = \alpha_{-}, K \in {\Bbb N}$,
then the empty word is in ${\mathcal D}(\sigma_-,\alpha_-)$
where
$\sigma_- = x_i$.
Otherwise let $j > i$ be given by 
$x_j \ne \alpha_-, x_{j+K} = \alpha_-, K \in {\Bbb N}$
and have
$x_{[i,j]} \in {\mathcal D}(\sigma_-,\alpha_-).$
The finiteness of
${\mathcal D}(\sigma_{-},\alpha_{-}), \, \sigma_{-}\in \Sigma_{-}(X)$
follows from condition $(c^-)$.
\end{proof}

Let $X \subset \Sigma^{\Bbb Z}$ be a subshift 
with a characteristic pair 
$((\alpha_{-})_{i\in {\Bbb Z}}, 
(\alpha_{+})_{i\in {\Bbb Z}})
$
of fixed points.
Let $x \in O_X$.
If for some $i_\circ \in {\Bbb Z}$,
\begin{align*}
x_i & = \alpha_-, \qquad i \le i_\circ,\\
x_i & = \alpha_+, \qquad i > i_\circ
\end{align*}
then set $c_X$ equal to the empty word.
Otherwise determine 
$i_-, i_+ \in {\Bbb Z}, i_- < i_+$, by
\begin{align*}
x_i & = \alpha_-, \qquad i < i_-,\\
x_{i_-}& \ne \alpha_-, \\
x_{i_+}&  \ne \alpha_+,\\
x_i & = \alpha_+, \qquad i > i_+,
\end{align*}
and set 
$c_X$ equal to the word 
$x_{[i_-,i_+]}$.

We also set 
\begin{align*}
\Omega^+(X) & = \{ d^+ \sigma_+ \mid \sigma_+ \in \Sigma_+(X), d^+ \in {\mathcal D}(\alpha_+,\sigma_+)\},\\
\Omega^-(X) & = \{ d_-^+ \sigma_+^- \mid \sigma_+^- \in \Sigma_+^-(X), d_-^+ \in {\mathcal D}(\alpha_-,\sigma_+^-)\}.
\end{align*}
We denote by 
$\OPR$ the set of 
$ d^+ \sigma_+ \in \Omega^+(X)$
such that there is a $D \in \Zp$
such that 
\begin{equation*}
\alpha_-^{k_-} c_X \alpha_+^{k_+ +D} d^+ \sigma_+ \in {\mathcal L}(X), \qquad k_-, k_+ \in {\Bbb N},
\end{equation*}
the smallest such $D$ to be denoted by $D(d^+\sigma_+)$.
We say that $X$ has reset if 
$\OPR \not= \emptyset.$
We set
$$
\OPC = \Omega^+(X) \backslash \OPR.
$$

We set
$$
\OMR =
\{c_X \alpha_+^{k_+ +D(d^+\sigma_+)} d^+ \sigma_+ \mid 
d^+ \sigma_+ \in \OPR,  k_+ \in \Bbb N\},
$$
and we say that $X$ satisfies the reset condition 
if
$\Omega^-(X) \backslash \OMR$ is a finite set.

For a strongly synchronizing subshift $X \subset \Sigma^{\Bbb Z}$
that has a synchronizing symbol denote by
${\mathcal C}(X)$ the set of admissible words of $X$
that begin with a synchronizing symbol,
that have no other synchronizing symbol
and that can be followed by a synchronizing symbol.
For $c \in {\mathcal C}(X)$
set $t(c)$ equal to the set of synchronizing symbols
that can follow $c$ and 
set $s(c)$ equal to the singleton set that contains the first symbol of $c$.
With the set of subsets of $\Sigma$ as index set and with a transition matrix
$A$ whose positive entries are given by
\begin{equation*}
A({\Sigma_\circ}, \{ \sigma\}) = 1, \qquad
\Sigma_\circ \in\bigl\{ \{ t(c) \mid c \in {\mathcal C}(X) \}\bigr\},\quad
 \sigma \in \Sigma_\circ, \quad
\sigma \in \synsyn, 
\end{equation*} 
${\mathcal C}(X)$
is a Markov code and
\begin{equation*}
X = scM({\mathcal C}(X)).
\end{equation*}

Given a strongly synchronizing subshift
$X \subset \Sigma^{\Bbb Z}$  
with a characteristic pair 
$((\alpha_{-})_{i\in {\Bbb Z}}, (\alpha_{+})_{i\in {\Bbb Z}})$
of fixed points
such that 
$\Sigma_-^+(X) = \emptyset$,
that satisfies the reset condition,
we set
\begin{equation}
\begin{split}
{\mathcal C}_-^{(X)}(\sigma_+^-) 
 = \{ \sigma_- d^- \alpha_-^{k_-} d_-^+ \mid
& \sigma_- \in \Sigma_-(X), d^- \in {\mathcal D}(\sigma_-,\alpha_-), 
\\ d_-^+\sigma_+^-  \in & \Omega^-(X)\backslash\Omega_{reset}^-(X), k_- \in {\Bbb N} \},
 \qquad
\sigma_+^- \in \Sigma_+^-(X), 
 \end{split}\label{eqn:41}
\end{equation}
\begin{align}
 {\mathcal C}_-^{(X)} 
& = \bigcup_{\sigma_+^- \in \Sigma_+^-(X)}
{\mathcal C}_-^{(X)}(\sigma_+^-), \label{eqn:42}  \\
t(c)
& = \{\sigma_+^- \in \Sigma_+^-(X) \mid
c \in {\mathcal C}_-^{(X)}(\sigma_+^-) \},\qquad c \in {\mathcal C}_-^{(X)}, \label{eqn:43}
\end{align}
and, given $M_-, M_+ \in \Zp$ and mappings 
\begin{align*}
\sigma_-d_- \longrightarrow & D_-(\sigma_- d^-)\in \Zp, \qquad 
\sigma_- \in \Sigma_-(X), d^-\in {\mathcal D}(\sigma_-,\alpha_-),\\
d^+ \sigma_+ \longrightarrow &  D_+(d^+ \sigma_+)\in \Zp, \qquad 
d^+ \sigma_+\in \OPR,
\end{align*} 
we set
\begin{gather}
\begin{split}
& {\mathcal C}_{reset}^{(X)}(D_-, M_-, M_+, D_+; \sigma_+) \\
= & \{ \sigma_- d^- \alpha_-^{k_-} c_X \alpha_+^{k_+ + D(d^+\sigma_+)} d^+ 
\mid \sigma_-\in \Sigma_-(X), d^- \in {\mathcal D}(\sigma_-,\alpha_-), d^+\sigma_+ \in \OPR,
\\ 
& k_{-}, k_{+} \in {\Bbb N},  
 D_-(\sigma_-d^-) + k_- + M_-  \ge M_+ + k_+ + D_+(d^+\sigma_+) \},\quad \sigma_+\in \Sigma_+(X), 
\end{split}\label{eqn:44} \\
{\mathcal C}_{reset}^{(X)}(D_-, M_-, M_+, D_+) 
 = \bigcup_{\sigma_+ \in \Sigma_+(X)}
{\mathcal C}_{reset}^{(X)}(D_-,M_-,M_+, D_+; \sigma_+), \label{eqn:45}  \\
\begin{split}
t(c)
& = \{\sigma_+ \in \Sigma_+(X) \mid
c \in {\mathcal C}_{reset}^{(X)}(D_-, M_-,M_+,D_+;\sigma_+) \},\\
& \qquad 
c \in {\mathcal C}_{reset}^{(X)}(D_-, M_-,M_+,D_+). 
\end{split}
\label{eqn:46}
\end{gather}
and given $J_-, J_+ \in \Zp$, and mappings 
\begin{align*}
\sigma_-d^- \longrightarrow & \Delta_-(\sigma_- d^-)\in \Zp, \qquad 
\sigma_- \in \Sigma_-(X), d^-\in {\mathcal D}(\sigma_-,\alpha_-),\\
d^+\sigma_+ \longrightarrow & \Delta_+(d^+ \sigma_+)\in \Zp, \qquad 
d^+\sigma_+ \in \OPC,
\end{align*} 
we set
\begin{gather}
\begin{split}
& {\mathcal C}_{counter}^{(X)}(\Delta_-, J_-, J_+, \Delta_+; \sigma^+) \\
= &  \{ \sigma_- d^- \alpha_-^{k_-} c_X \alpha_+^{k_+} d^+ \mid 
 \sigma_- \in \Sigma_-(X), d^- \in {\mathcal D}(\sigma_-,\alpha_-), \\
 & d^+\sigma_+ \in \OPC,
  k_{-}, k_{+} \in {\Bbb N}, \\
 & ( \Delta_-(\sigma_-, d^-) + k_- + J_- ) \cap (J_+ + k_+ + \Delta_+(d^+\sigma_+)) \ne \emptyset \},
\quad \sigma_+ \in \Sigma_+(X)
\end{split} \label{eqn:47}\\
{\mathcal C}_{counter}^{(X)}(\Delta_-,J_-,J_+, \Delta_+) 
 = \bigcup_{\sigma_+ \in \Sigma_+(X)}
{\mathcal C}_{counter}^{(X)}(\Delta_-, J_-, J_+,\Delta_+; \sigma_+), \label{eqn:48}  \\
\begin{split}
t(c) = \{\sigma_+ \in \Sigma_+(X) \mid
& c \in {\mathcal C}_{counter}^{(X)}(\Delta_-, J_-, J_+,\Delta_+, \sigma_+) \},\\
& c \in {\mathcal C}_{counter}^{(X)}(\Delta_-,J_-,J_+, \Delta_+). \label{eqn:49}
\end{split}
\end{gather}
By (3.1-9) there is defined a Markov code
$$
{\mathcal C}_-^{(X)} \cup
{\mathcal C}_{reset}^{(X)}(D_-, M_-,M_+,D_+)\cup
{\mathcal C}_{counter}^{(X)}(\Delta_-,J_-,J_+, \Delta_+).
$$

We define a standard one-counter shift as 
a strongly synchronizing subshift $X \subset {\Sigma}^{\Bbb Z}$
that has a characteristic pair 
of fixed points,
such that
$\Sigma_-^+(X)$
is empty,
such that 
$X$ satisfies the reset condition,
and such that
 there exist $ I \in \Zp$ and parameters
$D_-, M_-,M_+,D_+$,
$\Delta_-, J_-, J_+, \Delta_+$ 
such that 
\begin{equation}
\begin{split}
&\{ c \in {\mathcal C}(X) \mid \ell(c) > I \} \\
=& \{ c \in {\mathcal C}_{-}^{(X)}  
\cup       {\mathcal C}_{reset}^{(X)} (D_-, M_-, M_+,D_+)
\cup       {\mathcal C}_{counter}^{(X)} (\Delta_-, J_-, J_+, \Delta_+) \mid
\ell(c) > I \},
\end{split}
\label{eqn:511}
\end{equation}
where the equality is understood as an equality of Markov codes.
If (3.10) holds then we say that 
$I, D_-, M_-, M_+, D_+, \Delta_-, J_-, J_+, \Delta_+$
are parameters of the standard 
one-counter shift $X \subset {\Sigma}^{\Bbb Z}$.
The parameters
$
\Delta_-, J_-, J_+, \Delta_+
$
can be missing, and in the case that $X$ has no reset
the parameters
$D_-, M_-, M_+, D_+$
are missing.

For a standard one-counter shift 
$X \subset \Sigma^{\Bbb Z}$
denote the smallest $I \in \Zp$
such that
(3.10)
holds by $I_X$,
and denote by
$D_-(X), M_-(X), M_+(X), D_+(X)$,
$\Delta_-(X)$, 
$J_-(X)$,
$J_+(X)$,
$\Delta_+(X)$
the uniquely determined parameters for $X$
that satisfy the normalization conditions
\begin{align*}
 \min{(M_-,M_+)}  
=& \min{(J_-,J_+)} 
= \min_{\sigma_-\in \Sigma_-(X), d^-\in {\mathcal D}(\sigma_-,d_-)}  D(\sigma_-d^-)\\
=& \min_{d^+\sigma_+ \in {\Omega^+(X)}}  D(d^+\sigma_+) 
=\min_{\sigma_- \in \Sigma_-(X), d^- \in {\mathcal D}(\sigma_-,\alpha_-)}
\Delta_-(\sigma_-d^-) \\ 
=&\min_{\{\sigma_+ \in \Sigma_+(X), d^+ \in {\mathcal D}(\alpha_+,\sigma_+)\mid
d^+\sigma_+ \in \OPC \}}\Delta_+(d^+\sigma_+) = 0.
\end{align*}
E. g. 
for
$sc({\mathcal C}^{(N)}_{reset}\cup{\mathcal C}^{(N)}_{counter})
$
the normalized parameters are given by
$I = M_- = M_+ = J_- = J_+ = 0$,
the range of 
$D_-$ and $D_+$ being $\{0 \}$,
and 
the range of 
$\Delta_-$ and $\Delta_+$ being $\{0 \}$.
We associate with 
a standard one-counter shift $X \subset \Sigma^{\Bbb Z}$
 the Markov code
\begin{align*}
{\mathcal C}^{(X)} =
& \{ c \in {\mathcal C}^{(X)} \mid \ell(c) \le I_X \} \\
& \bigcup \{ c \in {\mathcal C}^{(X)}_{-} \cup 
   {\mathcal C}^{(X)}_{reset}(D_-(X), M_-(X), M_+(X),D_+(X))\\
&\quad \cup
   {\mathcal C}^{(X)}_{counter}(\Delta_-(X), J_-(X), J_+(X),\Delta_+(X))
   \mid \ell(c) > I_X \},
\end{align*}
and find that $X$
has a distinguished presentation as the Markov coded system of ${\mathcal C}^{(X)}$.

There is a development 
that can be considered, at least partially,
as the converse.
One takes as a starting point a finite alphabet $\Sigma$,
a proper subset $\Sigma_{synchro}$ of $\Sigma$ and symbols 
$\alpha_-, \alpha_+ \not\in \Sigma$.
One also has to provide for 
some $I \in \Zp$ 
a Markov code all of whose words 
begin with a symbol in $\Sigma_{synchro}$, with the 
remaining symbols in $\Sigma \backslash \Sigma_{synschro}$
and that have length less than or equal to to $I$,
and one has to provide the other components,
$\Sigma_-, \Sigma_+,$ 
$D_-, M_-, M_+, D_+$,
$\Delta_-, J_-, J_+, \Delta_+ $
that are needed for the construction of a Markov code 
${\mathcal C}$
according to rules that imitate the content of
(3.1-9).
An additional requirement is that 
the symbols in $\Sigma_{synchro}$
are the only synchronizing symbols in $scM({\mathcal C})$
for which there is a test.
One arrives in this way at a standard one-counter shift
$scM({\mathcal C})$ such that
${\mathcal C}(scM({\mathcal C})) = {\mathcal C},$
and one observes a perfect reciprocity
between a class of Markov codes 
and the class of standard one-counter shifts.
\medskip

{\bf 3 b. Behavior under topological conjugacy}


In this subsection
we assume that we are given a strongly synchronizing subshift
$X \subset \Sigma^{\Bbb Z}$
with a characteristic pair 
$(({\alpha}_-)_{i \in {\Bbb Z}},
({\alpha}_+)_{i \in {\Bbb Z}})
$
of fixed points 
and a subshift 
$X \subset \Sigma^{\Bbb Z}$ 
together with a topological conjugacy
$\tilde{\varphi}:\tilde{X} \longrightarrow X$
that is given by a 
one-block map
$\tilde{\varPhi}: \tilde{\Sigma} \longrightarrow \Sigma$
such that 
\begin{equation}
\tilde{\varPhi}^{-1}(\synsyn) \subset \tilde{\Sigma}_{synchro}(\tilde{X}).
\end{equation}
We introduce notation that we use in this situation.
We set
$
(\tilde{\alpha}_-)_{i \in {\Bbb Z}} = \tilde{\varphi}^{-1}(
(\alpha_-)_{i \in {\Bbb Z}})
$,
$
(\tilde{\alpha}_+)_{i \in {\Bbb Z}} = \tilde{\varphi}^{-1}(
(\alpha_+)_{i \in {\Bbb Z}}).
$
 $[-L,L]$ will denote a coding window of
$\tilde{\varphi}^{-1}$
and $\varPhi$ will be a block map 
$\varPhi: {\mathcal L}_{2L + 1}(X) \longrightarrow \tilde{\Sigma}$
that gives $\tilde{\varphi}^{-1}$.
$Q \in {\Bbb N}$ will be chosen such that 
for a synchronizing word $a$ of $X$ and for 
$a^- \in \Gamma_Q^-(a), 
a^+ \in \Gamma_Q^+(a)$
the word 
$a^- a a^+$ contains a synchronizing symbol.

For
$\tilde{\sigma}_- \in \tilde{\Sigma}_-(\tilde{X}),
\tilde{b}^-\in \Gamma_{Q + L}^-(\tilde{\sigma}_-),$
we can set by Lemma 2.1 and by (3.11)
$$
\tilde{\varPhi}(\tilde{b}^-\tilde{\sigma}_-) = 
b^-(\tilde{b}^-\tilde{\sigma}_-) \sigma_-(\tilde{b}^-\tilde{\sigma}_-) a^-(\tilde{b}^-\tilde{\sigma}_-),
$$
where the words 
$b^-(\tilde{b}^-\tilde{\sigma}_-) $
and 
$a^-(\tilde{b}^-\tilde{\sigma}_-) $
and the symbol
$\sigma_-(\tilde{b}^-\tilde{\sigma}_-) $
are uniquely determined by
$\tilde{b}^-\tilde{\sigma}_-$
 under the condition that 
$\sigma_-(\tilde{b}^-\tilde{\sigma}_-) $
is synchronizing and that 
$a^-(\tilde{b}^-\tilde{\sigma}_-) $
does not contain a synchronizing symbol.
We set 
$$
I_-(\tilde{b}^-\tilde{\sigma}_-) = \ell(a^-(\tilde{b}^-\tilde{\sigma}_-)).
$$
Denoting by 
$d^-(\tilde{b}^-\tilde{\sigma}_-\tilde{d}^-)$
the longest prefix of
$a^-(\tilde{b}^-\tilde{\sigma}_-)\tilde{\varPhi}(\tilde{d}^-)\alpha_-$
that is in ${\mathcal D}(\sigma_-,\alpha_-)$,
we have a mapping
$$
\Psi_{\tilde{b}^-\tilde{\sigma}_-}: \tilde{d}^- \longrightarrow d^-(\tilde{b}^-\tilde{\sigma}_-\tilde{d}^-),
\qquad
\tilde{d}^- \in {\mathcal D}(\tilde{\sigma}_-,\tilde{\alpha}_).
$$
Denote by 
${\mathcal D}_{\tilde{b}^-\tilde{\sigma}_-}(\alpha_-)$
the set of 
$d^- \in 
{\mathcal D}(\sigma_-(\tilde{b}^-\tilde{\sigma}_-),\alpha_-)
$
such that 
the prefix of length $Q + L +1$
of the word
$
b^-(\tilde{b}^-\tilde{\sigma}_-)\sigma_-(\tilde{b}^-\tilde{\sigma}_-)d^-\alpha_-^{Q+L}
$
is equal to 
$\tilde{\varPhi}(\tilde{b}^-\tilde{\sigma}_-)$
and such that the prefix of length
$Q + 2L + 1$ 
of the word
$
\varPhi(b^-(\tilde{b}^-\tilde{\sigma}_-)\sigma_-(\tilde{b}^-\tilde{\sigma}_-)d^-\alpha_-^{Q+2L})
$
is a suffix of 
$\tilde{b}^-\tilde{\sigma}_-$.
We use corresponding symbols with a time symmetric meaning.

We define 
$ H_-, H_+ \in \Zp$ 
by
\begin{equation}
\tilde{\varPhi}(c_{\tilde{X}}) 
= \alpha_-^{H_-} c_X \alpha_+^{H_+}.     
\end{equation}
\begin{lemma}
For 
$\tilde{\sigma}_- \in \tilde{\Sigma}_-(\tilde{X}),
\tilde{b}^- \in \Gamma_{Q+L}^-(\tilde{\sigma}_-),
$
the mapping
$\Psi_{\tilde{b}^-\tilde{\sigma}_-}$
is a bijection
of
${\mathcal D}(\tilde{\sigma}_-, \tilde{\alpha}_-)$
onto
${\mathcal D}_{\tilde{b}^-\tilde{\sigma}_-}(\alpha_-)$.
\end{lemma}
\begin{proof}
By construction
$$
\Psi_{\tilde{b}^-\tilde{\sigma}_-}
({\mathcal D}(\tilde{\sigma}_-, \tilde{\alpha}_-)) \subset 
{\mathcal D}_{\tilde{b}^-\tilde{\sigma}_-}(\alpha_-),
$$
and one confirms that the inverse of 
$
\Psi_{\tilde{b}^-\tilde{\sigma}_-}
$
is given by the mapping that assigns to a 
$d^- \in {\mathcal D}_{\tilde{b}^-\tilde{\sigma}_-}(\alpha_-)
$
the word that is obtained by removing the prefix of length $Q+1$ from the 
longest prefix of the word 
$
\varPhi(b^-(\tilde{b}^-\tilde{\sigma}_-)\sigma_-(\tilde{b}^-\tilde{\sigma}_-)d^-\alpha_-^{2L+1})
$
that does not end in $\tilde{\alpha}_-$.
\end{proof}
\begin{lemma}
Let
$\Sigma_-^+(X) = \emptyset$.
Then also 
$\tilde{\Sigma}_-^+(\tilde{X}) = \emptyset$.
\end{lemma}
\begin{proof}
If there were a 
$
\tilde{\sigma}^+_-\in \tilde{\Sigma}_-^+(\tilde{X})
$
and a 
$
\tilde{d}_+^- \in {\mathcal D}(\tilde{\sigma}_-^+, \tilde{\alpha}_+),
$
then one would have for a 
$\tilde{b}^- \in \Gamma^-_{Q+L}(\tilde{\sigma}_-^+)$
that
$$
\Psi_{\tilde{b}^-\tilde{\sigma}_-}(\tilde{d}_+^-)
 \in {\mathcal D}(\sigma_-^+(\tilde{b}^-\tilde{\sigma}_-), \alpha_+).
$$
\end{proof}
\begin{lemma}
For 
$\tilde{d}^+\tilde{\sigma}_+ \in \Omega_{reset}^+(\tilde{X}),
\tilde{b}^+ \in \Gamma_{L+Q}^+(\tilde{\sigma}_+),
$
one has
$$
d^+(\tilde{d}^+\tilde{\sigma}_+ \tilde{b}^+)
\sigma_+(\tilde{\sigma}_+\tilde{b}^+) \in \Omega^+_{reset}(X).
$$ 
\end{lemma}
\begin{proof}
One has 
$$
D_+(d^+(\tilde{d}^+\tilde{\sigma}_+\tilde{b}^+) \sigma_+(\tilde{\sigma}_+\tilde{b}^+)) 
\le
D_+(\tilde{d}^+\tilde{\sigma}_+) + I_+(\tilde{\sigma}_+\tilde{b}^+)
+
\ell(\tilde{d}^+) - 
\ell(d^+(\tilde{d}^+\tilde{\sigma}_+\tilde{b}^+)).
$$
\end{proof}
\begin{lemma}
Let 
$\tilde{d}_-^+\tilde{\sigma}_+^- \in \Omega^-(\tilde{X}),
\tilde{b}^+ \in \Gamma_{L+Q}^+(\tilde{\sigma}_+^-),
$
and let
\begin{align}
d^+\sigma_+(\tilde{\sigma}_+\tilde{b}^+) & \in \Omega^+_{reset}(X) \\
\intertext{and}
k_+ & > 2L 
\end{align}
be such that 
\begin{equation}
d_-^+(\tilde{d}_-^+\tilde{\sigma}_+^- \tilde{b}^+)
= c_X \alpha_+^{k_+ +D_+(d^+\sigma_+(\tilde{\sigma}_+\tilde{b}^+))}d^+.
\end{equation} 
Then
$$
\tilde{d}_-^+\tilde{\sigma}_+^- \in \Omega^-_{reset}(\tilde{X}).
$$ 
\end{lemma}
\begin{proof}
By (3.13)
$$
\alpha_-^{L + l_- + H_-}c_X 
\alpha_+^{l_+ + k_+ + D_+(d^+\sigma_+(\tilde{\sigma}_+\tilde{b}^+))}
d^+ \sigma_+(\tilde{\sigma}_+\tilde{b}^+) \in {\mathcal L}(X),
\quad l_-, l_+ \in {\Bbb N},
$$
and by (3.14) and (3.15) 
there is a $\tilde{k}_+ \in {\Bbb N}$
such that the word
$$
\varPhi(
\alpha_-^{L + l_- + H_-}c_X 
\alpha_+^{l_+ + k_+ + D_+(d^+\sigma_+(\tilde{\sigma}_+\tilde{b}^+))}
d^+ \sigma_+(\tilde{\sigma}_+\tilde{b}^+) b^+(\tilde{\sigma}_+\tilde{b}^+))
$$
contains the word
$$
\tilde{\alpha}_-^{l_-} \tilde{c}_X \tilde{\alpha}_+^{l_+ + \tilde{k}_+}
\tilde{d}^+_-\tilde{\sigma}^-_+
$$
as a subword for
$l_-, l_+ \in {\Bbb N}.$
\end{proof}
\begin{lemma}
For 
$\tilde{d}^+\tilde{\sigma}_+ \in \Omega^+(\tilde{X}),
\tilde{b}^+ \in \Gamma_{L+Q}^+(\tilde{\sigma}_+),
$
one has
$$
\tilde{d}^+\tilde{\sigma}_+ \in \Omega^+_{reset}(\tilde{X})
$$ 
if and only if
$$
d^+(\tilde{d}^+\tilde{\sigma}_+ \tilde{b}^+)
\sigma_+(\tilde{\sigma}_+\tilde{b}^+) \in \Omega^+_{reset}(X).
$$ 
\end{lemma}
\begin{proof}
This follows from Lemma 3.5 and Lemma 3.6.
\end{proof}
\begin{lemma}
Let $X$ satisfy the reset condition.
Then $\tilde{X}$ also satisfies the reset condition.
\end{lemma}
\begin{proof}
It follows from  Lemma 3.6
that there is a bound on the length of the words in
$
\Omega^-(\tilde{X}) \backslash \Omega^-_{reset}(\tilde{X}).
$
\end{proof}
We note that the converse of Lemma 3.8 also holds.
\begin{prop}
$\tilde{X}$ has reset if and only if $X$ has reset.
\end{prop}
\begin{proof}
Let
$d^+\sigma_+ \in \Omega^+(X)$.
To obtain 
$\tilde{\sigma}_+ \in \Sigma_+(\tilde{X})$
and
$\tilde{b}^+ \in \Gamma_{L+Q}^+(\tilde{\sigma}_+)$
such that 
$d^+\sigma_+ \in {\mathcal D}_{\tilde{\sigma}_+\tilde{b}^+}(\alpha_+)$,
let
$a^+ \in \Gamma_{L+Q}^+(\sigma_+)$
and let 
$\tilde{\sigma}_+\tilde{b}^+$
equal to the first subword of length $Q + L + 1$
of
$$
\varPhi(\alpha_+^{2L+1} d^+ \sigma_+ a^+)
$$
that begins with a synchronizing symbol.
Apply Lemma 3.3 and Lemma 3.7.
\end{proof}
\begin{lemma}
Let $X$ be a standard one-counter shift.
Then $\tilde{X}$ is also a standard one-counter shift.
\end{lemma}
\begin{proof}
By Lemma 3.4 
$\Sigma_-^+(\tilde{X})$ is empty and by Lemma 3.8
$\tilde{X}$ satisfies the reset condition.
Let
$I, D_-, M_-, M_+, D_+$,
$\Delta_-, J_-, J_+, \Delta_+$
be parameters for $X$.
Let
\begin{equation}
\begin{split}
\tilde{I} >
& I + 2Q + + 6L + M_- + M_+ + J_- + J_+ + \ell(C_X)   \\
& + 2\max \{ \ell(\sigma_- d^-) \mid \sigma_- \in \Sigma_-(X), d^-\in {\mathcal D}(\sigma_-,\alpha_-) \} \\
& + 2\max \{ \ell(d^+ \sigma_+) \mid d^+ \sigma_+\in \Omega^+(X) \} \\
& +  \max \{ \ell(d^+ \sigma_+^-) \mid d^+\sigma_+ \in \Omega^-(X) \backslash \Omega^-_{reset}(X) \} \\
& +  \max  \cup_{\sigma_- \in \Sigma_-(X), d^-\in {\mathcal D}(\sigma_-,\alpha_-)} \Delta_-(\sigma_-d^-) \\
& +  \max  \cup_{ d^+\sigma_+ \in \Omega^+_{counter}(X)} \Delta_+(d^+\sigma_+), 
\end{split}
\end{equation}
\begin{align}
\tilde{M}_- & = M_- + H_-, \qquad \tilde{M}_+ = M_+ + H_+, \\
\tilde{J}_- & = J_- + H_-, \qquad \tilde{J}_+ = J_+ + H_+, 
\end{align}
\begin{gather}
\begin{split}
{\tilde D}_-(\tilde{\sigma}_-\tilde{d}_-) 
& = \max_{\tilde{b}^- \in \tilde{\Gamma}^-_{Q+2L}(\tilde{\sigma}_-)}
\{
 D_-(\sigma_-(\tilde{b}^-\tilde{\sigma}_-) d^-(\tilde{b}^-\tilde{\sigma}_-\tilde{d}^-))  \\
& \qquad  -\ell(d^-(\tilde{b}^-\tilde{\sigma}_-\tilde{d}^-)) + \ell(\tilde{d}^-) + I_-(\tilde{b}^-\tilde{\sigma}_-)\},\\
& \qquad \qquad
\tilde{\sigma}_- \in \Sigma_-(\tilde{X}), 
\tilde{d}^- \in {\mathcal D}(\tilde{\sigma}_-,\tilde{\alpha}_-),
\end{split} \\
\begin{split}
\tilde{D}_+(\tilde{d}^+\tilde{\sigma}_+) 
& = \min_{\tilde{b}^+ \in \Gamma^+_{L+Q}(\tilde{\sigma}_+)}
\{ 
I_+(\tilde{\sigma}_+\tilde{b}^+)+ \ell(\tilde{d}^+)
-\ell(d^+(\tilde{d}^+\tilde{\sigma}_+\tilde{b}^+))  \\
& \qquad + 
D_+(d^+(\tilde{d}^+\tilde{\sigma}_+\tilde{b}^+)\sigma_+(\tilde{\sigma}_+\tilde{b}^+)) 
 \},\\
& \qquad \qquad 
\tilde{d}^+\tilde{\sigma}_+ \in \Omega^+_{reset}(\tilde{X}),
\end{split} 
\end{gather}
and
\begin{gather}
\begin{split}
\tilde{\Delta}_-(\tilde{\sigma}_-\tilde{d}^-)
& = \bigcup_{\tilde{b}^- \in \Gamma^-_{Q +L}(\tilde{\sigma}_-)}
 \Delta_-(\sigma_-(\tilde{b}^-\tilde{\sigma}_-)
    d^-(\tilde{b}^-\tilde{\sigma}_-\tilde{d}^-)) + \ell({\tilde{d}^-}) \\
& \qquad 
 -\ell(d^-(\tilde{b}^-\tilde{\sigma}_-\tilde{d}^-)) 
 + I_-(\tilde{b}^-\tilde{\sigma}_-)),\\
& \qquad \qquad
 \tilde{\sigma}_- \in \tilde{\Sigma}_-(\tilde{X}),
 \tilde{d}^- \in {\mathcal D}(\tilde{\sigma}_-, \tilde{\alpha}_-), 
\end{split} \\
\begin{split}
\tilde{\Delta}_+(\tilde{d}^+\tilde{\sigma}_+)
& = \bigcup_{\tilde{b}^+ \in \Gamma^+_{L+Q}(\tilde{\sigma}_+)}
 I_+(\tilde{\sigma}_+\tilde{b}^+)-\ell(d^+(\tilde{d}^+\tilde{\sigma}_+\tilde{b}^+)) + \ell(\tilde{d}^+) \\
& \qquad + \Delta_+(d^+(\tilde{d}^+\tilde{\sigma}_+\tilde{b}^+)\sigma_+(\tilde{\sigma}_+\tilde{b}^+)),\\
& \qquad \qquad
 \tilde{\sigma}_+ \tilde{d}^+ \in \Omega^+_{counter}(\tilde{X}).
\end{split}
\end{gather}
We prove that 
\begin{equation*}
\begin{split}
& \{ \tilde{c} \in {\mathcal C}(\tilde{X}) \mid \ell(\tilde{c})=\tilde{I} \}  \\
\subset \quad 
&  
{\mathcal C}^{(\tilde{X})}_{-} 
\cup 
 {\mathcal C}^{(\tilde{X})}_{reset}(\tilde{D}_-,\tilde{M}_-,\tilde{M}_+, \tilde{D}_+) 
\cup
{\mathcal C}^{(\tilde{X})}_{counter}(\tilde{\Delta}_-,\tilde{J}_-,\tilde{J}_+, \tilde{\Delta}_+).
\end{split}
\end{equation*}
Given a word $\tilde{c} \in {\mathcal C}(\tilde{X})$
of length $\tilde{I}$,
let
$\tilde{\sigma}_-$ be the first symbol of $\tilde{c}$ and let 
$\tilde{\sigma}_+ \in t(\tilde{c})$.
Also let 
$$
\tilde{b}^-\in \Gamma^-_{Q + L}(\tilde{\sigma}_-),\qquad
\tilde{b}^+\in \Gamma^+_{L + Q}(\tilde{\sigma}_+),
$$
and let a word $c \in {\mathcal L}(X)$
be given by 
$$
\tilde{\varPhi}(\tilde{b}^- \tilde{c}\tilde{\sigma}_+ \tilde{b}^+)
=b^- ( \tilde{b}^- \tilde{\sigma}_-) c \sigma_+(\tilde{\sigma}_+ \tilde{b}^+)
 b^+(\tilde{\sigma}_+ \tilde{b}^+).
 $$
By (3.16)
\begin{equation*}
\begin{split}
c \in &  
{\mathcal C}^{(X)}_{-}(\sigma_+(\tilde{\sigma}_+ \tilde{b}^+)) 
\cup 
 {\mathcal C}^{(X)}_{reset}(D_-,M_-,M_+, D_+;\sigma_+(\tilde{\sigma}_+ \tilde{b}^+))\\ 
& \cup
{\mathcal C}^{(X)}_{counter}(\Delta_-,J_-,J_+, \Delta_+;\sigma_+(\tilde{\sigma}_+ \tilde{b}^+)).
\end{split}
\end{equation*}
In the case that 
\begin{equation}
c \in 
{\mathcal C}^{(X)}_{-}(\sigma_+(\tilde{\sigma}_+ \tilde{b}^+)), 
\end{equation}
one has 
$
\sigma_+(\tilde{\sigma}_+ \tilde{b}^+) \in \Sigma_+^-(X),
$
and there are
$$
d^- \in {\mathcal D}(\sigma_-(\tilde{b}^-\tilde{\sigma}_-),\alpha_-),\qquad
d^+_- \in {\mathcal D}(\alpha_-,\sigma_+(\tilde{\sigma}_+ \tilde{b}^+)),
$$
and $k_-\in {\Bbb N}$ such that
\begin{equation*}
d^+_-\sigma_+(\tilde{\sigma}_+ \tilde{b}^+) \in \Omega^-(X)\backslash\Omega_{reset}^-(X),
\qquad
c=\sigma_-(\tilde{b}^-\tilde{\sigma}_-) d^- \alpha_-^{k_-} d_-^+.  
\end{equation*}
By (3.16)
$k_- > 2 L$,
and it is seen from the action of $\varPhi$ 
that one has,  setting
$$
\tilde{d}^- = \Psi^{-1}_{\tilde{b}^-\tilde{\sigma}_-}(d^-),
\qquad
\tilde{d}_-^+ = \Psi^{-1}_{\tilde{\sigma}_+\tilde{b}^+}(d_-^+),
$$
and 
\begin{equation}
\tilde{k}_- = k_- - I_-(\tilde{b}^-\tilde{\sigma}_-) + \ell(d^-) - 
\ell(\tilde{d}^-) - H_-,
\end{equation}
that 
\begin{equation*}
\tilde{c} 
= \tilde{\sigma}_- \tilde{d}^- 
\tilde{\alpha}_-^{\tilde{k}_-}
\tilde{d}^+_-. 
\end{equation*}
Here 
$$
\tilde{d}_-^+ \tilde{\sigma}_+ 
\in \Omega^-(\tilde{X})\backslash\Omega^-_{reset}(\tilde{X}),
$$
for otherwise one would have by Lemma 3.5 a contradiction to (3.23).
This means that 
$$
\tilde{c} \in {\mathcal C}_-^{(\tilde{X})}(\tilde{\sigma}_+).
$$

In the case that 
$$
c \in
 {\mathcal C}^{(X)}_{reset}(D_-,M_-,M_+, D_+;\sigma_+(\tilde{\sigma}_+ \tilde{b}^+)),
$$
there are
$$
d^- \in {\mathcal D}(\sigma_-(\tilde{b}^-\tilde{\sigma}_-),\alpha_-),\qquad
d^+_- \in {\mathcal D}(\alpha_+,\sigma_+(\tilde{\sigma}_+ \tilde{b}^+)),
$$
and $k_-, k_+ \in {\Bbb N}$ such that
$$
d^+ \sigma_+(\tilde{\sigma}_+ \tilde{b}^+) \in \Omega^+_{reset}(X),
\qquad
c = \sigma_-(\tilde{b}^-\tilde{\sigma}_-) d^- \alpha_-^{k_-}c_X
\alpha_+^{k_+} d^+,
$$
\begin{equation}
D_-(\sigma_-(\tilde{b}^-\tilde{\sigma}_-)d^-) + k_- + M_- 
\ge 
M_+ + k_+ + D_+(d^+\sigma_+(\tilde{\sigma}_+ \tilde{b}^+)).  
\end{equation}
Set again
\begin{align}
\tilde{d}^- & = \Psi^{-1}_{\tilde{b}^-\tilde{\sigma}_-}(d^-),\\
\intertext{and also set}
\tilde{d}^+ & = \Psi^{-1}_{\tilde{\sigma}^+\tilde{b}^+}(d^+),
\end{align}
$$
\tilde{d}_-^+ = \Psi^{-1}_{\tilde{\sigma}_+\tilde{b}^+}(c_X \alpha_+^{k_+}d^+).$$
If here
\begin{equation}
\tilde{d}_-^+ \tilde{\sigma}_+ 
\in \Omega^-(\tilde{X})\backslash\Omega^-_{reset}(\tilde{X}),
\end{equation}
then by Lemma 3.6
$k_+ \le 2 L$,
and then by (3.16)
$k_- > 2L$,
and it is seen from the action of $\varPhi$ that, 
with $\tilde{k}_-$ given by the expression (3.24), 
\begin{equation*}
\tilde{c} 
= \tilde{\sigma}_- \tilde{d}^- 
\tilde{\alpha}_-^{\tilde{k}_-}
\tilde{d}^+_-. 
\end{equation*}
By (3.28) this means that
$$
\tilde{c} \in {\mathcal C}_-^{(\tilde{X})}(\tilde{\sigma}_+).
$$

If here
$$
\tilde{d}_-^+ \tilde{\sigma}_+ 
\in \Omega^-_{reset}(\tilde{X}),
$$
one has by (3.16) and (3.25)
that $k_- > 2L$ 
and it is seen from the action of $\varPhi$ that,
with $\tilde{k}_-$ given by the expression (3.24),
and with
\begin{equation}
\tilde{k}_+ = k_+ - H_- - \ell(\tilde{d}^+) + \ell(d^+) - I_+(\tilde{\sigma}_+\tilde{b}^+),
\end{equation}
that one has then
\begin{equation*}
\tilde{c} 
= 
\tilde{\sigma}_- \tilde{d}^- 
\tilde{\alpha}_-^{\tilde{k}_-} c_{\tilde{X}}
\tilde{\alpha}_+^{\tilde{k}_+}
\tilde{d}^+.
\end{equation*}
By (3.25), (3.19) and (3.20)
\begin{equation*}
\tilde{D}_-(\tilde{\sigma}_- \tilde{d}^-) + \tilde{k}_- + \tilde{M}_- 
\ge
\tilde{M}_+ + \tilde{k}_+ + \tilde{D}_+(\tilde{d}^+\tilde{\sigma}_+),
\end{equation*}
and this means that 
$$
\tilde{c} \in {\mathcal C}^{(\tilde{X})}_{reset}(\tilde{D}_-,\tilde{M}_-,\tilde{M}_+, \tilde{D}_+;\tilde{\sigma}_+). 
$$

In case that 
$$
c \in
{\mathcal C}^{(X)}_{counter}(\Delta_-, J_-, J_+, \Delta_+;\sigma_+),
$$
there are
$$
d^- \in {\mathcal D}(\sigma_-(\tilde{b}^-\tilde{\sigma}_-),\alpha_-),\qquad
d^+ \in {\mathcal D}(\alpha_+,\sigma_+(\tilde{\sigma}_+ \tilde{b}^+)),
$$
and
$$
D_- \in \Delta_-(\sigma_-(\tilde{b}^-\tilde{\sigma}_-)d^-),\qquad
D_+ \in \Delta_+(d^+\sigma_+(\tilde{\sigma}_+ \tilde{b}^+)),
$$
and $k_-, k_+ \in {\Bbb N}$ 
such that
\begin{equation}
d^+\sigma_+(\tilde{\sigma}_+ \tilde{b}^+) \in \OPC,
\end{equation}
\begin{equation*}
 c=  \sigma_-(\tilde{b}^-\tilde{\sigma}_-) d^- \alpha_-^{k_-} c_X {\alpha_+}^{k_+}d^+,
\end{equation*}
\begin{equation}
 D_- + k_- + J_- = J_+ + k_+ + D_+.
\end{equation}
By (3.16) and (3.31)
$k_-, k_+ > 2 L$,
and with 
$\tilde{d}^-,\tilde{d}^+,\tilde{k}_-,\tilde{k}_+$
given by the expressions (3.26),(3.27), (3.24), (3.29) and with 
\begin{align}
\tilde{D}_-& = D_- + \ell(\tilde{d}^-) -\ell(d^-) + I_-(\tilde{b}^-\tilde{\sigma}_-),\\
\tilde{D}_+& = I_+(\tilde{\sigma}_+ \tilde{b}^+) + \ell(\tilde{d}^+) -\ell(d^+) + D_+,
\end{align}
it is seen from the action of $\varPhi$ that 
\begin{equation*}
\tilde{c} 
= 
\tilde{\sigma}_- \tilde{d}^- 
\tilde{\alpha}_-^{\tilde{k}_-} c_{\tilde{X}}
\tilde{\alpha}_+^{\tilde{k}_+}
\tilde{d}^+.
\end{equation*}
By (3.32) and (3.33)
\begin{equation*}
\tilde{D}_- + \tilde{k}_- + \tilde{J}_- 
=
\tilde{J}_+ + \tilde{k}_+ + \tilde{D}_+,
\end{equation*}
and by Lemma 3.7 and by (3.30)
this means that
$$
\tilde{c} \in
{\mathcal C}^{(\tilde{X})}_{counter}(\tilde{\Delta}_-, \tilde{J}_-, \tilde{J}_+, \tilde{\Delta}_+;\tilde{\sigma}_+).
$$

We prove that 
\begin{equation}
\{ \tilde{c} \in {\mathcal C}_-^{(\tilde{X})} \mid \ell(\tilde{c}) = \tilde{I} \} 
\subset 
{\mathcal C}(\tilde{X}).
\end{equation}
For 
$\tilde{\sigma}_-\in \tilde{\Sigma}_+^-(\tilde{X}),$
and for a word
$ 
 \tilde{c} \in {\mathcal C}_-^{(\tilde{X})}(\tilde{\sigma}_+^-)
 $
 of length $\tilde{I}$,
 with the first symbol $\tilde{\sigma}_-$,
there are
$$
\tilde{d}^- \in {\mathcal D}(\tilde{\sigma}_-,\tilde{\alpha}_-),\qquad
\tilde{d}^+_- \in {\mathcal D}(\tilde{\alpha}_+,\tilde{\sigma}_+^-),
$$
and $\tilde{k}_- \in {\Bbb N}$ such that
\begin{equation}
 \tilde{d}_-^+ \tilde{\sigma}_+^- 
\in \Omega^-(\tilde{X})\backslash\Omega^-_{reset}(\tilde{X}),
\end{equation}
\begin{equation*}
 \tilde{c} 
= \tilde{\sigma}_- \tilde{d}^- 
\tilde{\alpha}_-^{\tilde{k}_-}
\tilde{d}^+_-. 
\end{equation*}
Let
$$
\tilde{b}^-\in \Gamma^-_{Q + L}(\tilde{\sigma}_-),\qquad
\tilde{b}^+\in \Gamma^+_{Q + L}(\tilde{\sigma}_+),
$$
and let a word $c \in {\mathcal L}(X)$
in the symbols of $\Sigma$ be given by 
\begin{equation}
\varPhi(\tilde{b}^- \tilde{c}\tilde{\sigma}_+ \tilde{b}^+)
=b^- ( \tilde{b}^- \tilde{\sigma}_-) c \sigma_+^-(\tilde{\sigma}_+^- \tilde{b}^+).
\end{equation}
From (3.36)
it is seen that  
there is a
 $k_-\in {\Bbb N}$ such that
\begin{equation*}
c=\sigma_-(\tilde{b}^-\tilde{\sigma}_-) d^-(\tilde{b}^-\tilde{\sigma}_-\tilde{d}^-) \alpha_-^{k_-} 
d_-^+(\tilde{d}^+_-\tilde{\sigma}_+^- \tilde{b}^+).  
\end{equation*}
If here
$$
d_-^+(\tilde{d}^+_-\tilde{\sigma}_+^- \tilde{b}^+)
\sigma_+^-(\tilde{\sigma}^-_+\tilde{b}^+)  \in \Omega^-(X)\backslash\Omega_{reset}^-(X),
$$
then by (3.16), $k_- > 2L$ 
and if here
$$
d_-^+(\tilde{d}^+_-\tilde{\sigma}_+^- \tilde{b}^+)
\sigma_+^-(\tilde{\sigma}^-_+\tilde{b}^+)
 \in \Omega_{reset}^-(X),
$$
then there are
$
d^+ \in {\mathcal D}(\alpha_+,\sigma_+^-(\tilde{\sigma}_+^- \tilde{b}^+))
$
and $k_+ \in {\Bbb N}$ such that 
$$
d^+\sigma_+^- \in \OPR,
\qquad
d_-^+(\tilde{d}^+_-\tilde{\sigma}_+^- \tilde{b}^+)
= c_X\alpha_+^{k_+}d^+.
$$
By Lemma 3.6 and by (3.35)
$k_+ \le 2L$,
and then by (3.16) 
$k_- >2L$, and also
\begin{equation*}
D_-(\sigma_-(\tilde{b}^-\tilde{\sigma}_-)d^-) + k_- + M_- 
\ge 
M_+ + k_+ + D_+(d^+\sigma_+(\tilde{\sigma}_+ \tilde{b}^+)),  
\end{equation*}
and therefore 
$$
c c_X\alpha_+^{k_+}d^+ \in 
{\mathcal C}^{(X)}_{reset}(D_-,M_-,M_+, D_+;\sigma_+(\tilde{\sigma}_+ \tilde{b}^+)).
$$
By (3.16)
then
$$
c c_X\alpha_+^{k_+}d^+ \in {\mathcal C}(X),
$$
and it is seen from the action of $\varPhi$
that the word $\tilde{c}$
is a subword of the word
$$
\varPhi(b^-(\tilde{b}^-\tilde{\sigma}_-)c c_X\alpha_+^{k_+}d^+\sigma_+(\tilde{\sigma}_+ \tilde{b}^+) 
b^+(\tilde{\sigma}_+ \tilde{b}^+))
 \in {\mathcal L}(X),
$$
and (3.34) is confirmed.

We prove that
\begin{equation}
\{ \tilde{c} \in
 {\mathcal C}^{(\tilde{X})}_{reset}(\tilde{D}_-,\tilde{M}_-,\tilde{M}_+, \tilde{D}_+)
 \mid \ell(\tilde{c} ) = \tilde{I} \} \subset {\mathcal C}(\tilde{X}).
\end{equation}
For $\tilde{\sigma}_+ \in \tilde{\Sigma}_+(\tilde{X})$
and for a word
$$
\tilde{c} \in
 {\mathcal C}^{(\tilde{X})}_{reset}(\tilde{D}_-,\tilde{M}_-,\tilde{M}_+, \tilde{D}_+;\tilde{\sigma}_+)
$$
of length $\tilde{I}$
with the first symbol
$\tilde{\sigma}_-$
 there are
$$
\tilde{d}^- \in {\mathcal D}(\tilde{\sigma}_-, \tilde{\alpha}_-),\qquad
\tilde{d}^+ \in {\mathcal D}(\tilde{\alpha}_+, \tilde{\sigma}_+),
$$
and $k_-, k_+ \in {\Bbb N}$ such that
\begin{align}
\tilde{d}_-^+ \tilde{\sigma}_+ & \in \Omega^-_{reset}(\tilde{X}),\\
\tilde{D}_-(\tilde{\sigma}_-\tilde{d}^-) + \tilde{k}_- + \tilde{M}_- 
& \ge 
\tilde{M}_+ + \tilde{k}_+ + \tilde{D}_+(\tilde{d}^+\tilde{\sigma}_+),
\end{align}
\begin{equation*}
\tilde{c} 
= \tilde{\sigma}_- \tilde{d}^- 
\tilde{\alpha}_-^{\tilde{k}_-}
c_{\tilde{X}}
\tilde{\alpha}_+^{\tilde{k}_+}
\tilde{d}^+. 
\end{equation*}
By (3.12), (3.17), (3.19) and (3.20) one can select
$$
\tilde{b}^-\in \Gamma^-_{Q + L}(\tilde{\sigma}_-),\qquad
\tilde{b}^+\in \Gamma^+_{Q + L}(\tilde{\sigma}_+)
$$
such that
\begin{align*}
\tilde{D}_-(\tilde{\sigma}_-\tilde{d}^-)
& = D_-(\sigma_-(\tilde{b}^-\tilde{\sigma}_-)d^-(\tilde{b}^-\tilde{\sigma}_-\tilde{d}^-))
- \ell(d^-(\tilde{b}^-\tilde{\sigma}_-\tilde{d}^-)) + \ell(\tilde{d}^-) + I_-(\tilde{b}^-\tilde{\sigma}_-),\\
\tilde{D}_+(\tilde{d}^+\tilde{\sigma}_+)
& = I_+(\tilde{\sigma}_+ \tilde{b}^+) + \ell(\tilde{d}^+) -\ell(d^+(\tilde{d}^+\tilde{\sigma}_+\tilde{b}^+)) +
 D_+(d^+(\tilde{d}^+\tilde{\sigma}_+\tilde{b}^+)\sigma_+(\tilde{\sigma}_+ \tilde{b}^+)),
\end{align*}
and such that
one has with
\begin{align}
k_- & = I_-(\tilde{b}^-\tilde{\sigma}_-) -\ell(d^-(\tilde{b}^-\tilde{\sigma}_-\tilde{d}^-)) + \ell(\tilde{d}^-)
+ \tilde{k}_- + H_-, \\
k_+ & = H_+ + k_+ +\ell(\tilde{d}^+) -\ell(d^+(\tilde{d}^+\tilde{\sigma}_+\tilde{b}^+)) +
 I_+(\tilde{\sigma}_+\tilde{b}^+),
\end{align}
that 
\begin{equation}
D_-(\sigma_-(\tilde{b}^-\tilde{\sigma}_-)d^-(\tilde{b}^-\tilde{\sigma}_-\tilde{d}^-))
+ k_- + M_- 
\ge 
M_+ + k_+ +
 D_+(d^+(\tilde{d}^+\tilde{\sigma}_+\tilde{b}^+)\sigma_+(\tilde{\sigma}_+ \tilde{b}^+)).
\end{equation}
By (3.38) and Lemma 3.7 and by (3.42)
it follows for the word $c$ in the symbols of $\Sigma$ that is given by 
\begin{equation*}
\tilde{\varPhi}(\tilde{b}^- \tilde{c}\tilde{\sigma}_+ \tilde{b}^+)
=
b^- ( \tilde{b}^- \tilde{\sigma}_-) c \sigma_+(\tilde{\sigma}_+ \tilde{b}^+)b^+(\tilde{\sigma}_+ \tilde{b}^+),
\end{equation*}
that
\begin{equation*}
c
=
\sigma_-(\tilde{b}^-\tilde{\sigma}_-) d^-(\tilde{b}^-\tilde{\sigma}_-\tilde{d}^-) \alpha_-^{k_-} 
 c_X \alpha_+^{k_+}d^+(\tilde{d}^+\tilde{\sigma}_+^- \tilde{b}^+)
\in {\mathcal C}^{(X)}_{reset}(\sigma_+(\tilde{\sigma}_+ \tilde{b}^+)).  
\end{equation*}
By (3.16) then 
$c \in {\mathcal C}(X)$.
By (3.16) and (3.42)
$k_-, k_+ > 2L$
and from the action of $\varPhi$ it is seen that the word 
$\tilde{c}$ is a subword of the word
\begin{equation*}
{\varPhi}(
b^- ( \tilde{b}^- \tilde{\sigma}_-) c \sigma_+(\tilde{\sigma}_+ \tilde{b}^+)b^+(\tilde{\sigma}_+ \tilde{b}^+)
) \in {\mathcal L}(X),
\end{equation*}
and
(3.37) is confirmed.

We prove that
\begin{equation}
\{ \tilde{c} \in
 {\mathcal C}^{(\tilde{X})}_{counter}(\tilde{\Delta}_-,\tilde{J}_-,\tilde{J}_+, \tilde{\Delta}_+)
 \mid \ell(\tilde{c} ) = \tilde{I} \} \subset {\mathcal C}(\tilde{X}).
\end{equation}
For $\tilde{\sigma}_+ \in \tilde{\Sigma}_+(\tilde{X})$
and a word
$$
\tilde{c} \in
 {\mathcal C}^{(\tilde{X})}_{counter}(\tilde{\Delta}_-,\tilde{J}_-,\tilde{J}_+, \tilde{\Delta}_+;\tilde{\sigma}_+)
$$
of length $\tilde{I}$
with a first symbol
$\tilde{\sigma}_-$
 there are
$$
\tilde{d}^- \in {\mathcal D}(\tilde{\sigma}_-, \tilde{\alpha}_-),\qquad
\tilde{d}^+ \in {\mathcal D}(\tilde{\alpha}_+, \tilde{\sigma}_+),
$$
and 
$$
\tilde{D}^- \in \tilde{\Delta}_-(\tilde{\sigma}_- \tilde{d}^-),\qquad
\tilde{D}^+ \in \tilde{\Delta}_+(\tilde{d}^+ \tilde{\sigma}_+),
$$
and
$\tilde{k}_-, \tilde{k}_+ \in {\Bbb N}$ such that
\begin{align}
\tilde{d}^+ \tilde{\sigma}_+ & \in \Omega^-_{counter}(\tilde{X}),\\
\tilde{D}_- + \tilde{k}_- + \tilde{J}_- 
& = 
\tilde{J}_+ + \tilde{k}_+ + \tilde{D}_+.
\end{align}
By (3.12), (3.18), (3.21) and (3.22) one can select
$$
\tilde{b}^-\in \Gamma^-_{Q + L}(\tilde{\sigma}_-),\qquad
\tilde{b}^+\in \Gamma^+_{Q + L}(\tilde{\sigma}_+),
$$
such that there are 
\begin{align*}
D_- \in
&  \Delta_-(\sigma_-(\tilde{b}^-\tilde{\sigma}_-)d^-(\tilde{b}^-\tilde{\sigma}_-\tilde{d}^-)),\\
D_+ \in 
&  \Delta_+(d^+(\tilde{d}^+\tilde{\sigma}_+\tilde{b}^+)\sigma_+(\tilde{\sigma}_+ \tilde{b}^+)),
\end{align*}
such that
\begin{align*}
\tilde{D}_-
& = D_- - \ell(d^-(\tilde{b}^-\tilde{\sigma}_-\tilde{d}^-)) + \ell(\tilde{d}^-) + I_-(\tilde{b}^-\tilde{\sigma}_-),\\
\tilde{D}_+
& = I_+(\tilde{\sigma}_+ \tilde{b}^+) + \ell(\tilde{d}^+) -\ell(d^+(\tilde{d}^+\tilde{\sigma}_+\tilde{b}^+)) +
 D_+.
\end{align*}
With $k_-, k_+ \in {\Bbb N}$
given by the expressions (3.40) and (3.41) then
\begin{equation}
D_- + k_- + J_- = J_+ + k_+ + D_+.
\end{equation}
By (3.44) and Lemma 3.7 and by (3.45)
it follows for  the word $c$ in the symbls of $\Sigma$ that is given by 
\begin{equation*}
\tilde{\varPhi}(\tilde{b}^- \tilde{c}\tilde{\sigma}_+ \tilde{b}^+)
=
b^- ( \tilde{b}^- \tilde{\sigma}_-) c \sigma_+(\tilde{\sigma}_+ \tilde{b}^+)b^+(\tilde{\sigma}_+ \tilde{b}^+),
\end{equation*}
that
\begin{equation*}
c
=
\sigma_-(\tilde{b}^-\tilde{\sigma}_-) d^-(\tilde{b}^-\tilde{\sigma}_-\tilde{d}^-) \alpha_-^{k_-} 
 c_X \alpha_+^{k_+}d^+(\tilde{d}^+\tilde{\sigma}_+^- \tilde{b}^+)
\in {\mathcal C}^{(X)}_{counter}(\sigma_+(\tilde{\sigma}_+ \tilde{b}^+)).  
\end{equation*}
By (3.16) then 
$c \in {\mathcal C}(X)$.
By (3.16) and (3.46)
$k_-, k_+ > 2L$
and from the action of $\varPhi$ it is seen that the word 
$\tilde{c}$ is a subword of the word
\begin{equation*}
{\varPhi}(
b^- ( \tilde{b}^- \tilde{\sigma}_-) c \sigma_+(\tilde{\sigma}_+ \tilde{b}^+)b^+(\tilde{\sigma}_+ \tilde{b}^+)
) \in {\mathcal L}(\tilde{X}),
\end{equation*}
and
(3.43) is confirmed.

We have shown that 
$
\tilde{I}, \tilde{D}_-,\tilde{M}_-,\tilde{M}_+, \tilde{D}_+,
\tilde{\Delta}_-, \tilde{J}_-, \tilde{J}_+, \tilde{\Delta}_+
$
are parameters for $\tilde{X}$.
 \end{proof}


\medskip

{\bf 3 c. Shifts of standard one-counter type}
    
One has a theorem that can be viewed as analogous to Theorem 1.1.
\begin{theorem}
Let $X \subset \Sigma^{\Bbb Z}$ be a subshift 
that is topologically conjugate to a standard one-counter shift.
Then there exists an  $n_\circ \in {\Bbb N}$
such that 
$X^{{\langle [1,n]\rangle}}$ is a standard one-counter shift,
\begin{equation*}
X^{{\langle [1,n]\rangle}} = scM({\mathcal C}^{(X^{{\langle [1,n]\rangle}})}), 
\qquad n \ge n_\circ.
\end{equation*}
\end{theorem}
\begin{proof}
Apply Lemma 2.2 and Lemma 3.10.
\end{proof}
One can view the class of standard one-counter shifts 
as extending the class of topological Markov shifts 
and one is then lead to introduce a class of subshifts of 
standard one-counter type as the class of subshifts
that have a  higher block system
that is a standard one-counter shift.
Theorem 3.11 is then equivalent to the statement 
that a subshift that is topologically conjugate to a subshift 
of standard one-counter type is itself of standard one-counter type.

\section{$\lambda$-graph systems and $C^*$-algebras}

Consider a $\lambda$-graph system ${\frak L} =(V,E,\lambda,\iota)$ 
over the alphabet $\Sigma$ 
with vertex set
$
V = \cup_{l \in \Zp} V_{l},
$
edge set
$
E = \cup_{l \in \Zp} E_{l,l+1},
$
labeling map 
$\lambda: E \rightarrow \Sigma$
and shift-like map 
$
\iota
$
that is given by surjective maps
$
\iota_{l,l+1}:V_{l+1} \rightarrow V_l,l \in \Zp.
$
A subset $\V$ of $V$ is called hereditary
if all $v \in V$ such that $\iota(v) \in \V$ are in $\V$,
and if $v \in \V$ then all initial vertices of all edges that have $v$ as a final vertex 
are also in $\V$.
A hereditary subset $\V$ is said to be proper if 
$\V \cap V_l \ne V_l$ for all $l \in {\Bbb N}$.

Let us denote by 
$\{v_1^l,\dots, v_{m(l)}^l\}$
the vertex set
$V_l$ at level $l$.
For
$
i=1,2,\dots,m(l),\ j=1,2,\dots,m(l+1), \ \alpha \in \Sigma,
$ 
we put
\begin{eqnarray*}
A_{l,l+1}(i,\alpha,j)
  &=&
{\begin{cases}
1 & \text{ if } \ s(e) = {v}_i^l, \lambda(e) = \alpha,
      t(e) = {v}_j^{l+1} \text{ for some }    e \in E_{l,l+1}, \\
0           & \text{ otherwise,}
\end{cases}} \\
I_{l,l+1}(i,j)
  &=&
{\begin{cases}
1 &  
    \text{ if } \ \iota_{l,l+1}({v}_j^{l+1}) = {v}_i^l, \\
0           & \text{ otherwise.}
\end{cases}} 
\end{eqnarray*}
The $C^*$-algebra 
${\mathcal O}_{\frak L}$ 
associated with ${\frak L}$ 
is the  universal    
$C^*$-algebra generated by partial isometries
$S_\alpha, \alpha \in \Sigma$
and projections $E_i^l, i=1,2,\dots,m(l),\ l \in \Zp$  
subject to the following  operator relations called $({\frak L})$:
\begin{eqnarray*}
& &\sum_{\beta \in \Sigma}  S_{\beta}S_{\beta}^*  =  1, \\
  \sum_{i=1}^{m(l)} E_i^l   &=&  1, \qquad 
 E_i^l   =  \sum_{j=1}^{m(l+1)}I_{l,l+1}(i,j)E_j^{l+1}, \\
& & S_\alpha S_\alpha^* E_i^l =E_i^l 
S_\alpha S_\alpha^*,  \\ 
S_{\alpha}^*E_i^l S_{\alpha}  &=&  
\sum_{j=1}^{m(l+1)} A_{l,l+1}(i,\alpha,j)E_j^{l+1},
\end{eqnarray*}
for
$
i=1,2,\dots,m(l),\l\in \Zp, 
 \alpha \in \Sigma
$
\cite{Ma2002a}.

For a subshift $X\subset \Sigma^{\Bbb Z}$
we recall the construction of its future $\lambda$-graph system
${}^X\!{\frak L}$.
The label set of 
${}^X\!{\frak L}$ is $\Sigma$
and its vertex set is
$$
V(X) = \cup_{l \in \Zp}V_l(X)
$$
where $V_0(X)$ contains the singleton set that contains the empty word,
and where
$$
V_l(X) = \{ \Gamma_l^+(x^-) \mid x^- \in X_{(-\infty,0]} \}, \qquad
l \in {\Bbb N}.
$$
All edges of $\PL$ leave a vertex in $\cup_{l \in {\Bbb N}} V_l(X)$,
and a vertex $ v \in \cup_{l \in {\Bbb N}} V_l(X)$
has an outgoing edge that carries the label $\sigma \in \Sigma$
if and only if $v$ contains a word that begins with $\sigma$
and the target vertex of this outgoing edge 
is equal to 
$\{ a \in \Gamma^+_{[1,l)} \mid \sigma a \in v \}, l \in {\Bbb N}$.
The mapping
$$
\iota : \cup_{l \in {\Bbb N}} V_l(X) \longrightarrow \cup_{l \in \Zp} V_l(X)
$$
deletes last symbols.

\begin{theorem}
Let $X \subset \Sigma^{\Bbb Z}$
be a standard one-counter shift
with a characteristic pair 
$(({\alpha}_-)_{i \in {\Bbb Z}},
({\alpha}_+)_{i \in {\Bbb Z}})
$
of fixed points.
Then
\begin{enumerate}\renewcommand{\labelenumi}{(\roman{enumi})}
\item $V(X)$ has a proper hereditary subset if and only if $X$ has no reset.
\item $\FL$ has a proper hereditary subset.
\end{enumerate}
\end{theorem}
\begin{proof}
(i) Let $\OPR \ne \emptyset$.
Let $I, D_-, M_-, M_+, D_+$ be parameters for $X$,
where
$I$ is chosen such that 
$scM(\{ c \in {\mathcal L}(X) \mid \ell(c) \le I \})$
is aperiodic and 
topologically transitive subshift of finite type with alphabet $\Sigma$.
Let $Q \in {\Bbb N}$ 
be such that for $\sigma, \sigma' \in \synsyn$
there exists for $ q > Q$ an admissible concatenation of words in
$\{ c \in {\mathcal L}(X) \mid \ell(c) \le I \}$
that begins with $\sigma$ and that can be followed by $\sigma'$.
With $D > I$ such that also 
$$
D > \ell(c_X) + M_- + M_+ + \ell(d^-) + D_-(\sigma_-d^-), \qquad
\sigma_- \in \Sigma_-, d^- \in {\mathcal D}(\sigma_-,\alpha_-),
$$ 
one has for $x^- \in X_{(-\infty,0]}$,
that
$\Gamma_D^+(x^-)$
contains a synchronizing symbol.
Let $x^- \in X_{(-\infty,0]}, l \in {\Bbb N}$.
One can choose a word $ a\in {\mathcal L}(X)$ of length less than $l + D$ such that
$$
\Gamma_l^+(x^-) =
\Gamma_l^+(a) 
$$
and
for $y^- \in X_{(-\infty,0]}$
one has that 
$\Gamma_{l +2D +Q}^+(y^-)$ contains a word with suffix $a$.
It follows that $V(X)$ has no proper hereditary subset.

In case that $\OPR = \emptyset$ one has
$\{ \alpha_+^l \} \in V_l(X), l \in {\Bbb N}$
and it follows that the set 
$\cup_{l \in {\Bbb N}} V_l(X) \backslash \{ \alpha_+^l \}$
is a proper hereditary subset of $V(X)$.

(ii)
Here $\cup_{l \in {\Bbb N}}(V_l(X) \backslash \{ \alpha_-^l \}) $
is a proper hereditary subset of $V(X)$.
\end{proof}

\begin{cor}
Let $X$ be a standard one-counter shift.
Then
\begin{enumerate}\renewcommand{\labelenumi}{(\roman{enumi})}
\item ${\mathcal O}_{\PL}$ is simple if and only if $X$ has reset.
\item ${\mathcal O}_{\FL}$ is not simple.
\end{enumerate}
\end{cor}
\begin{proof}
There exists a bijective correspondence between hereditary subsets of the vertex set $V$
and ideals in the  $C^*$-algebra ${\mathcal O}_{\frak L}$
(\cite{Ma2002a}, \cite{Ma2006a}).
\end{proof}

For the notion of flow equivalence see \cite{BF, Fr, PS, Th}.
For a subshift 
$Y \subset \Sigma^{\Bbb Z}$
and for 
$
\sigma \in \Sigma,\
 \sigma' \not\in \Sigma,
$
we say that the subshift
$
Y' \subset (\Sigma \cup \{ \sigma' \})^{\Bbb Z}
$
is obtained from the subshift $Y$
by replacing in $Y$
$\sigma$ by $\sigma \sigma'$
if
for every admissible word $a'$ of $Y'$
there is an admissible word $a$ of $Y$
such that
$a'$ can be obtained by replacing in $a$ the symbol $\sigma$
by the word $\sigma \sigma'$
and then,
 if necessarily,
still removing the first symbol
or the last symbol or both.
We say then also that $Y$ 
is obtained from $Y'$ by replacing in $Y'$ 
the word $\sigma \sigma'$ 
by the symbol $\sigma$.

Subshifts
$X \subset \Sigma^{\Bbb Z}$
and
$\widetilde{X} \subset \widetilde{\Sigma}^{\Bbb Z}$
are flow equivalent if there exists a chain of subshifts
\begin{equation*}
Y[q] \subset \Sigma[q]^{\Bbb Z}, \quad 1 \le q \le Q, \ Q \in {\Bbb N}, 
\qquad Y[1] = X, \quad Y[Q] = \widetilde{X},
\end{equation*}
such that
$Y[q]$ is topologically conjugate to 
$Y[q+1]$
or
$Y[q+1]$
is obtained from 
$Y[q]$ 
by replacing in $Y[q]$ 
a symbol $\sigma$
by the word
$\sigma \sigma'$
or
$Y[q]$ 
is obtained from 
$Y[q+1]$ 
by replacing in $Y[q+1]$ 
a symbol $\sigma$
by the word
$\sigma \sigma'$.
We remark at this point that the definition of a standard  
one-counter shift can be given a more general formulation in which the  
characteristic pair of fixed points are replaced by a pair of periodic points. 
In this  
way one arrives at a class of subshifts that is closed under  
flow equivalence.
\begin{cor}
A standard one-counter shift with reset 
is not flow equivalent to its inverse.
\end{cor}
\begin{proof}
The ideal structure of 
the $C^*$-algebra ${\mathcal O}_{\PL}$ is an invariant of flow equivalence
\cite{Ma2001b}. 
Apply Theorem 4.1. 
\end{proof}

\section{K-groups}
We will compute the K-groups and the Bowen-Franks groups of 
the one-counter shift 
\begin{equation*}
\sccnrt 
= sc(\{ a_n \alpha_+^m \alpha_-^k \mid 1 \le n\le N, \ m,k \in {\Bbb N}, m \le k \}),
\end{equation*}
that is of the  future $\lambda$-graph system of $\sccnrt$ or,
equivalently, of the  past $\lambda$-graph system of $\sccnr$.
The set up that we choose is for the future $\lambda$-graph system of $\sccnrt$.  
  Let $(\M,I) =(\M_{l,l+1}, I_{l,l+1})_{l\in \Zp}$ 
be the  symbolic matrix system of 
$\sccnrt$ (the  future $\lambda$-graph system of 
$\sccnrt$).
Let 
$(M,I) = (M_{l,l+1}, I_{l,l+1})_{l\in \Zp}$ 
 be its nonnegative matrix system.
 The entries of the nonnegative matrix $M_{l,l+1}$ 
count the number of symbols of the corresponding entries of 
 $\M_{l,l+1}$.
We denote by $m(l)$ the row size of $\M_{l,l+1}$,
so that the both matrices $M_{l,l+1}$ and $ I_{l,l+1}$
are $m(l) \times m(l+1)$ matrices. 
They satisfy the following relations
$$
 I_{l,l+1}M_{l+1,l+2} = M_{l,l+1} I_{l+1,l+2}, \qquad l \in \Zp. 
$$
We denote by
$
   \bar{I}^t_{l,l+1}, l \in \Zp
$
the homomorphism from
$
{\Bbb Z}^{m(l)} / (M_{l-1,l}^t - I_{l-1,l}^t){{\Bbb Z}^{m(l-1)}}
$
to
$
{\Bbb Z}^{m(l+1)} / (M_{l,l+1}^t - I_{l,l+1}^t){{\Bbb Z}^{m(l)}}
$
induced by 
${I}^t_{l,l+1}.$
Then as in \cite{Ma1999c}
\begin{align}
K_0( \sccnrt) 
& 
= \underset{l}{\varinjlim} 
\{ {\Bbb Z}^{m(l+1)} / (M_{l,l+1}^t - I_{l,l+1}^t){{\Bbb Z}^{m(l)}}, 
   \bar{I}^t_{l,l+1} \}, \\
K_1(\sccnrt)
&
=\underset{l}{\varinjlim}\{ \Ker (M_{l,l+1}^t - I_{l,l+1}^t)
             \text{ in }  {\Bbb Z}^{m(l)}, I^t_{l,l+1} \}. 
\end{align}
Let
${\Bbb Z}_{I}$ 
be the group of the projective limit 
$\underset{l}{\varinjlim} 
\{ {\Bbb Z}^{m(l)}, {I}_{l,l+1} \}.
$
The sequence 
$M_{l,l+1} - I_{l,l+1}, l \in \Zp$
acts on it as an endomorphism, denoted by $M - I.$
The Bowen-Franks groups $BF^i(\sccnrt ), i=0,1,$ are defined by
\begin{align*}
BF^0(\sccnrt) & = {\Bbb Z}_{I} / (M-I){\Bbb Z}_{I}, \\
BF^1(\sccnrt) & = \Ker(M-I) \quad \text{ in }\quad {\Bbb Z}_{I}.
\end{align*}
 We denote the symbols $\alpha_+, \alpha_-$ 
in the subshifts $\sccnrt$
now by $b,c$ respectively. 
  For $l \in {\Bbb N}$,
 consider the following subsets $\{ F_i^l \}_{i=1, \dots, 2l+2}$
 of the right one-sided shift $\sccnrt_{[1,\infty)}$.
\begin{align*}
F_1^l = & \{ {(x_n)}_{n\in {\Bbb N}} \in \sccnrt_{[1,\infty)} \mid 
x_1 = b, x_2 = \cdots = x_{l+2} = c,\\
&\qquad \qquad \qquad \qquad \qquad \qquad  
x_{l+3} = a_i \text{ for some } 1\le i \le N \},\\   
F_2^l = & \{ {(x_n)}_{n\in {\Bbb N}} \in \sccnrt_{[1,\infty)} \mid 
x_1 = b, x_2 = \cdots = x_{l+1} = c,\\
&\qquad \qquad \qquad \qquad \qquad \qquad   
x_{l+2} = a_i \text{ for some } 1\le i \le N \},\\   
  & \vdots \\
F_{l+1}^l = & \{ {(x_n)}_{n\in {\Bbb N}} \in \sccnrt_{[1,\infty)} \mid 
x_1 = b, x_2 = c, x_{3} = a_i \text{ for some } 1\le i \le N \},\\   
F_{l+2}^l = & \{ {(x_n)}_{n\in {\Bbb N}} \in \sccnrt_{[1,\infty)} \mid 
x_1 =  a_i \text{ for some } 1\le i \le N \},\\   
F_{l+3}^l = & \{ {(x_n)}_{n\in {\Bbb N}} \in \sccnrt_{[1,\infty)} \mid 
x_1 = c, x_2 = a_i \text{ for some } 1\le i \le N \},\\   
F_{l+4}^l = & \{ {(x_n)}_{n\in {\Bbb N}} \in \sccnrt_{[1,\infty)} \mid 
x_1 = x_2 = c, x_3 = a_i \text{ for some } 1\le i \le N \},\\   
 & \vdots \\
F_{2l+2}^l = & \{ {(x_n)}_{n\in {\Bbb N}} \in \sccnrt_{[1,\infty)} \mid 
x_1 = \cdots x_l = c, x_{l+1}= a_i \text{ for some } 1\le i \le N \}.  
\end{align*}
The sets
$\{ F_i^l \}_{i=1,\dots,2l+2}$ 
are the $l$-past equivalence classes of $ \sccnrt$.
Put
$m(l) = 2l+2.$ 
Let $v_i^l, i=1,\dots,m(l)$ be the vertex set $V_l$ of the canonical
$\lambda$-graph system ${\frak L}^{\sccnrt}$for the subshift $\sccnrt$.
The vertex $v_i^l$ is considered to be the class $[F_i^l]$ of $F_i^l$.
For a symbol $\gamma$, if $\gamma F_j^{l+1} $ is contained in $F_i^l$, 
then a labeled edge labeled $\gamma$ from the vertex $v_i^l$ to the vertex 
$v_j^{l+1}$
is defined in the $\lambda$-graph system.
Hence
there are labeled edges labeled $a_n,n=1,\dots,N$ from $v_{l+2}^l$ to $v_j^{l+1}$ for $j=1,2,\dots, l+2$. 
There are labeled edges labeled $b$ from $v_{i}^l$ to $v_{2l+4-i}^{l+1}$
and to $v_i^{l+1}$ for $i=1,2,\dots, l+1$.
There are labeled edges labeled $c$ from $v_i^l$  to $v_i^{l+1}$ for $i = l+3,l+4,\dots, 2l+2$,
and from $v_{2l+2}^l $ to $v_{2l+3}^{l+1}$ and to $v_{2l+4}^{l+1}.$

If $F_j^{l+1} $ is contained in $F_i^l$,
the $\iota$-map is defined by $\iota(v_j^{l+1}) = v_i^l$.
Hence we have
\begin{equation*}
\iota(v_j^{l+1})
=
\begin{cases}
v_1^l & \text{ if } j=1,\\
v_{j-1}^l & \text{ if } j=2,3,\dots, 2l+3,\\
v_{2l+2}^l & \text{ if } j=2l + 4.
\end{cases}
\end{equation*}

We will consider the symbolic matrix system $\M_{l,l+1}, I_{l,l+1}$
on the ordered bases
$
F_1^l,  \cdots, F_{m(l)}^l.
$
For $i=1,\dots, m(l)$ and $j= 1,\dots,m(l+1)$, we have
\begin{align*}
\M_{l,l+1}(i,j) & =
{\begin{cases}
a_1 + \cdots + a_N & \text{ if } i=l+2, \, j= 1,2,\cdots, l+2,\\
b & \text{ if } 1 \le i= j \le l+1, \\
b & \text{ if } i + j = 2l +5, \, 1 \le i \le l+1, \\
c & \text{ if } l+3 \le i= j \le 2l+2, \\
c & \text{ if } i = 2l + 2, j= 2l +3, 2l+4, \\
0 & \text{ otherwise, } 
\end{cases}}\\
I_{l,l+1}(i,j)& =
{\begin{cases}
1 & \text{ if } i=j = 1, \\
1 & \text{ if } 2 \le j= i+1 \le 2l+3, \\
1 & \text{ if } i=2l+2,  j= 2l+4, \\
0 & \text{ otherwise. } 
\end{cases}}
\end{align*}
Hence we have 
$$
M_{l,l+1}^t(i,j) =
\begin{cases}
N & \text{ if } j=l+2, \, i= 1,2,\cdots, l+2,\\
1 & \text{ if } 1 \le i= j \le l+1, \\
1 & \text{ if } i + j = 2l +5, \, 1 \le j \le l+1, \\
1 & \text{ if } l+3 \le i= j \le 2l+2, \\
1 & \text{ if } j = 2l + 2, i= 2l +3, 2l+4, \\
0 & \text{ otherwise, } 
\end{cases}
$$  
so that
$$
M_{l,l+1}^t(i,j) - I_{l,l+1}^t(i,j)  =
\begin{cases}
N & \text{ if } j=l+2, \, i= 1,2,\cdots, l+2,\\
1 & \text{ if } 2 \le i= j \le 2l+2, \, i\ne l+2,\\
1 & \text{ if } i + j = 2l +5, \, 1 \le j \le l+1, \\
-1 & \text{ if } 2 \le i = j+1\le 2l + 2, \\
0 & \text{ otherwise, } 
\end{cases}
$$  

$$
\setcounter{MaxMatrixCols}{16}
M_{l,l+1}^t -I_{l,l+1}^t
=
\begin{bmatrix}
            0&                         \hdotsfor{5}     & 0   & N   & 0    & \hdotsfor{5}             \\
           -1& 1& 0    &               \hdotsfor{3}     & 0   & N   & 0    & \hdotsfor{5}             \\
            0&-1& 1    & 0     &       \hdotsfor{2}     & 0   & N   & 0    & \hdotsfor{5}             \\
\hdotsfor{1} & 0&-1    & 1     & 0    &\hdotsfor{1}     & 0   & N   & 0    & \hdotsfor{5}             \\
\hdotsfor{2}    & \cdot& \cdot &\cdot &\cdot            &\cdot&\cdot& \cdot& \hdotsfor{5}             \\
\hdotsfor{3}           & \cdot &\cdot &\cdot            &\cdot&\cdot& \cdot& \hdotsfor{5}             \\
\hdotsfor{4}                   & 0    & -1              & 1   & N   & 0    & \hdotsfor{5}             \\
\hdotsfor{5}                          & 0               & -1  & N   & 0    & \hdotsfor{5}             \\
\hdotsfor{6}                                            & 0   &-1   & 1    & 0   &              \hdotsfor{4} \\
\hdotsfor{5}                                      & 0   & 1   & 0   &-1    & 1   & 0    &       \hdotsfor{3} \\ 
\hdotsfor{4}                               & 0    & 1   & 0   & 0   & 0    &-1   & 1    & 0    &\hdotsfor{2} \\
\hdotsfor{3}           & \cdot & \cdot&\cdot      & \hdotsfor{3}    &\cdot &\cdot&\cdot &\cdot &\cdots       \\ 
\hdotsfor{2}    & \cdot& \cdot & \cdot&\hdotsfor{5}                        &\cdot&\cdot &\cdot &\cdot        \\
\hdotsfor{1} & 0& 1    &0      & \hdotsfor{7}                                    &0     & -1   &1            \\
            0& 1& 0    &         \hdotsfor{10}                                                 &0            \\
            1& 0&                \hdotsfor{11}                                                 &0 
 \end{bmatrix}. 
$$
We see that
\begin{lemma}
$
\Ker(M_{l,l+1}^t - I_{l,l+1}^t) 
= 0
\quad
\text{ for }
\quad 
2 \le l \in {\Bbb N}.
$
\end{lemma}
Thus we have by (5.2),
\begin{lemma}
$
K_1(\sccnrt)\cong 0. 
$
\end{lemma}

We will next compute $K_0(\sccnrt)$.
Set for $i=1,\dots, 2l+4, \, j=1,\dots, 2l +2$
$$
B_{l,l+1}(i,j)  =
\begin{cases}
N & \text{ if } j=l+2, \, i= 1,\\
1 & \text{ if } 2 \le i= j \le 2l+2, \, i\ne l+2,\\
1 & \text{ if } (i, j) = (2l +4, 1), (2l+3, 2), \\
-1 & \text{ if } 2 \le i = j+1\le l + 2, \\
0 & \text{ otherwise. } 
\end{cases}
$$  
That is 
$$
\setcounter{MaxMatrixCols}{16}
B_{l,l+1}
=
\begin{bmatrix}
            0&                         \hdotsfor{5}     & 0   & N   & 0    & \hdotsfor{5}             \\
           -1& 1& 0    &               \hdotsfor{3}     & 0   & 0   & 0    & \hdotsfor{5}             \\
            0&-1& 1    & 0     &       \hdotsfor{2}     & 0   & 0   & 0    & \hdotsfor{5}             \\
\hdotsfor{1} & 0&-1    & 1     & 0    &\hdotsfor{1}     & 0   & 0   & 0    & \hdotsfor{5}             \\
\hdotsfor{2}    & \cdot& \cdot &\cdot &\cdot            &\cdot&\cdot& \cdot& \hdotsfor{5}             \\
\hdotsfor{3}           & \cdot &\cdot &\cdot            &\cdot&\cdot& \cdot& \hdotsfor{5}             \\
\hdotsfor{4}                   & 0    & -1              & 1   & 0   & 0    & \hdotsfor{5}             \\
\hdotsfor{5}                          & 0               & -1  & 0   & 0    & \hdotsfor{5}             \\
\hdotsfor{6}                                            & 0   & 0   & 1    & 0   &              \hdotsfor{4} \\
\hdotsfor{5}                                      & 0   & 0   & 0   & 0    & 1   & 0    &       \hdotsfor{3} \\ 
\hdotsfor{4}                               & 0    & 0   & 0   & 0   & 0    & 0   & 1    & 0    &\hdotsfor{2} \\
\hdotsfor{3}           & \cdot & \cdot&\cdot      & \hdotsfor{3}    &\cdot &\cdot&\cdot &\cdot &\cdots       \\ 
\hdotsfor{2}    & \cdot& \cdot & \cdot&\hdotsfor{5}                        &\cdot&\cdot &\cdot &\cdot        \\
\hdotsfor{1} & 0& 0    &0      & \hdotsfor{7}                                    &0     &  0   &1            \\
            0& 1& 0    &         \hdotsfor{10}                                                 &0            \\
            1& 0&                \hdotsfor{11}                                                 &0 
 \end{bmatrix}. 
$$
Let $P_l$ 
be the $(2 l +2) \times (2 l+2)$ matrix defined by setting for $i,j= 1,\dots, 2l+2$,
$$
P_l(i,j)  =
\begin{cases}
1 & \text{ if } i= j,\\
-1 & \text{ if } j=1, \, i= 2, \dots, l+1, \\
0 & \text{ otherwise. } 
\end{cases}
$$  
We know that 
\begin{equation}
P_{l+1} (M_{l,l+1}^t - I_{l,l+1}^t) {\Bbb Z}^{2l+2} = B_{l,l+1}{\Bbb Z}^{2l+2}. 
\end{equation}
Denote by 
$\bar{P}_{l+1}$ the induced isomorphism from
${\Bbb Z}^{2l+4} /(M_{l,l+1}^t - I_{l,l+1}^t) {\Bbb Z}^{2l+2}$
onto
${\Bbb Z}^{2l+4} /B_{l,l+1}{\Bbb Z}^{2l+2}$.
Let
$J_{l,l+1}$ be the 
$(2 l +4) \times (2 l+2)$ matrix defined by setting for $i=1,\dots, 2l+4, \, j= 1,\dots, 2l+2$,
$$
J_{l,l+1}(i,j)  =
\begin{cases}
1 & \text{ if } i= j=1,\\
1 & \text{ if } i=j-1, \, i= 2, \dots, 2l+3, \\
1 & \text{ if } i=2l+4, \, j= 2l+2, \\
0 & \text{ otherwise. } 
\end{cases}
$$  
Denote by 
$\bar{J}_{l,l+1}$ the induced homomorphism from
${\Bbb Z}^{2l+2} /B_{l-1,l}{\Bbb Z}^{2l}$
into
${\Bbb Z}^{2l+4} /B_{l,l+1}{\Bbb Z}^{2l+2}$.
\begin{lemma}
The diagram 
$$
\begin{CD}
{\Bbb Z}^{2l+2}/(M_{l-1,l}^t - I_{l-1,l}^t) {\Bbb Z}^{2l}
@>\bar{I}^t_{l,l+1}>> 
{\Bbb Z}^{2l+4}/(M_{l,l+1}^t - I_{l,l+1}^t) {\Bbb Z}^{2l+2}
 \\
@V \bar{P}_l VV @V \bar{P}_{l+1} VV \\
{\Bbb Z}^{2l+2}/ B_{l-1,l}{\Bbb Z}^{2l}
@>\bar{J}_{l,l+1}>> 
{\Bbb Z}^{2l+4}/ B_{l,l+1}{\Bbb Z}^{2l+2}
\end{CD}
$$
is commutative.
\end{lemma}

For an integer $n$, 
we denote by $q(n) \in {\Bbb Z}$ 
the quotient of $n$ by $N$
and 
by $r(n) \in \{0,1,\dots, N-1\}$
its residue such as
$n = q(n) N + r(n)$.
The following lemma is straightforward.
\begin{lemma}
Fix $l=2,3,\dots $.
For
$z =
\begin{bmatrix}
z_1  \\
\vdots   \\
z_{2l+4}
\end{bmatrix} 
\in {\Bbb Z}^{2l+4},
$
put inductively
\begin{align*}
x_1       & = z_{2l+4}, \\
x_2       & = z_{2l+3}, \\
x_{k}     & = z_{k} \qquad \text{ for } k=l+3,l+4,\dots, 2l+2,\\
x_{l+2}   & = q(z_{1}),\\
x_{l+1}   & = -z_{l+2},\\
x_{l}     & = -z_{l+1} -z_{l+2},\\
x_{l-k}   & = - z_{l-k+1} - z_{l-k+2} - \cdots - z_{l+2}, \qquad \text{ for } k=1,2,\dots, l-3.
\end{align*}   
Set
\begin{align*}
r_{l,l+1}(z) & = r(z_1) \in \{ 0,1,\dots, N-1 \},\\
\varphi_{l,l+1}(z) & =  z_2 - z_{2l+3} + z_{2l+4},\\
\psi_{l,l+1}(z) & = z_3 + z_4+ z_5 + \cdots + z_{l+2} +z_{2l+3}.
\end{align*}
Then we have
\begin{equation*}
\begin{bmatrix}
z_1  \\
\vdots   \\
z_{2l+4}
\end{bmatrix} 
=
B_{l,l+1}
\begin{bmatrix}
x_1  \\
\vdots   \\
x_{2l+2}
\end{bmatrix} 
+
\begin{bmatrix}
r_{l,l+1}(z)  \\
\varphi_{l,l+1}(z) \\
\psi_{l,l+1}(z)\\
0 \\
\vdots   \\
0
\end{bmatrix}.
\end{equation*}
\end{lemma}
The following lemma is also direct.
\begin{lemma}
For
$z= [z_i]_{i=1}^{2l+4} \in {\Bbb Z}^{2l+4}$,
one has 
$$
r_{l,l+1}(z) = 0 \text{ in } \{0,1,\dots,N-1 \} 
\quad
\text{ and }
\quad
\varphi_{l,l+1}(z) = \psi_{l,l+1}(z)=0 \text{ in } {\Bbb Z} 
$$
if and only if
there exists
$y= [y_i]_{i=1}^{2l+2}\in {\Bbb Z}^{2l+2}$ 
such that 
$ z =B_{l,l+1}(y)$.
\end{lemma}
\begin{lemma}
The map
$
\xi_{l+1} : 
[z_i]_{i=1}^{2l+4} \in {\Bbb Z}^{2l+4} \longrightarrow
(r_{l,l+1}(z), \varphi_{l,l+1}(z), \psi_{l,l+1}(z)) 
\in \{0,1,\dots,N-1 \}\oplus{\Bbb Z}\oplus{\Bbb Z} 
$
induces an isomorphism from
$
{\Bbb Z}^{2l+4}/B_{l,l+1}{\Bbb Z}^{2l+2}
$
onto
$ 
{\Bbb Z}/N{\Bbb Z} \oplus {\Bbb Z}\oplus{\Bbb Z}.
$
\end{lemma}
\begin{proof}
It suffices to show the surjectivity of $\xi_{l+1}$.
For 
$(g,m,k) \in \{0,1,\dots,N-1 \}\oplus{\Bbb Z}\oplus{\Bbb Z}$,
put
$z =[g,m,k,0,\dots,0]^t \in {\Bbb Z}^{2l+4}$.
One then sees that
$$
r_{l,l+1}(z) = g, \qquad \varphi_{l,l+1}(z) =m, \qquad  \psi_{l,l+1}(z) =k.
$$
\end{proof}
We denote by 
$\bar{\xi}_{l+1}$
the above isomorphism
from
$
{\Bbb Z}^{2l+4}/B_{l,l+1}{\Bbb Z}^{2l+2}
$
onto
$ 
{\Bbb Z}/N{\Bbb Z} \oplus {\Bbb Z}\oplus{\Bbb Z}
$
induced by 
${\xi}_{l+1}$.
\begin{lemma}
The diagram 
$$
\begin{CD}
{\Bbb Z}^{2l+2}/B_{l-1,l}{\Bbb Z}^{2l}
@>\bar{J}_{l,l+1}>> 
{\Bbb Z}^{2l+4}/B_{l,l+1}{\Bbb Z}^{2l+2}
 \\
@V \bar{\xi}_l VV @V \bar{\xi}_{l+1} VV \\
{\Bbb Z}/N{\Bbb Z} \oplus {\Bbb Z}\oplus{\Bbb Z} 
@>L >> 
{\Bbb Z}/N{\Bbb Z} \oplus {\Bbb Z}\oplus{\Bbb Z}
\end{CD}
$$
is commutative,
where
$L
=
\begin{bmatrix}
1 & 0 & 0 \\
0 & 0 & 0 \\
0 & 1 & 1 
\end{bmatrix}.
$
\end{lemma}
\begin{proof}
For
$z = [z_i ]_{i=1}^{2l+2} \in {\Bbb Z}^{2l+2}$,
it is direct to see that
\begin{align*}
r_{l,l+1}(J_{l,l+1}(z)) & = r_{l-1,l}(z),  \qquad 
\varphi_{l,l+1}(J_{l,l+1}(z)) = 0,  \\ 
\psi_{l,l+1}(J_{l,l+1}(z)) & = \varphi_{l-1,l}(z) + \psi_{l-1,l}(z).
\end{align*}
\end{proof}

Therefore we  conclude
\begin{lemma}
 $K_0(\sccnrt ) \cong {\Bbb Z}/N{\Bbb Z} \oplus {\Bbb Z}$.
\end{lemma}
\begin{proof}
By Lemma 5.1, 
it follows that
\begin{align*}
K_0(sc(\sccnrt ) 
= & \varinjlim 
\{ {\Bbb Z}^{2l + 4} / B_{l,l+1}{\Bbb Z}^{2l+2}, \overline{I^t}_{l,l+1} \} \\
= & \varinjlim \{ {\Bbb Z}/N{\Bbb Z} \oplus {\Bbb Z}\oplus{\Bbb Z}, L \} \\
\cong & {\Bbb Z}/N{\Bbb Z} \oplus {\Bbb Z}.
\end{align*}
\end{proof}

As the torsion free part of 
$
K_0(\sccnrt)
$
is not isomorphic to 
$
K_1(\sccnrt)
$,
these types of K-groups cannot appear in those of sofic systems.

We next  compute   the Bowen-Franks groups 
$BF^0(\sccnrt)$ 
and 
$BF^1(\sccnrt).$
As in \cite[Theorem 9.6]{Ma1999c},  
one sees the following formulae of short exact sequences of the 
universal coefficient type theorem:
\begin{align*}
0& \rightarrow
   \Ext_{\Bbb Z}^1(K_{i}(\sccnrt),{\Bbb Z}) \\
 & \rightarrow
   BF^{i}(\sccnrt) \\
&  \rightarrow
   \Hom_{\Bbb Z}(K_{i+1}(\sccnrt),{\Bbb Z})
\rightarrow
0.
\end{align*}
The sequences split unnaturally. 
\begin{lemma}
$
BF^0(\sccnrt) \cong {\Bbb Z}/ N {\Bbb Z},
\quad
BF^1(\sccnrt) \cong {\Bbb Z}^2.
$
\end{lemma} 
\begin{proof}
Since for a finitely generated abelian group $G$, 
$\Hom_{\Bbb Z}(G,\Bbb Z)$ is the torsion free part of $G$
and
$\Ext_{\Bbb Z}^1(G,\Bbb Z)$
is the torsion part of $G$,
one gets the desired assertions by 
Lemma 5.8.
\end{proof} 
As the torsion free part of 
$
BF^0(\sccnrt)
$
is not isomorphic to 
$
BF^1(\sccnrt)
$,
these types of Bowen-Franks groups cannot appear in those of sofic systems.
We restate Lemma 5.2, Lemma 5.8 and Lemma 5.9 as 
\begin{theorem}
\begin{align*}
K_0(\sccnrt)  \cong & {\Bbb Z}/N {\Bbb Z} \oplus {\Bbb Z}, \qquad 
K_1(\sccnrt)  \cong  0,\\
BF^0(\sccnrt) \cong & {\Bbb Z}/N{\Bbb Z}, \quad \qquad
BF^1(\sccnrt) \cong {\Bbb Z}^2.
\end{align*}
\end{theorem}

We will next compute the K-groups for 
$\sccnr$.
The computation is completely similar to the above one as in the following way.
We can take the $l$-past equivalence classes of $\sccnr $
as the similar ones to the $\sccnrt$.
Let $(\M,I) = (\M_{l,l+1},I_{l,l+1})_{l\in \Zp}$
be the symbolic matrix system for $\sccnr$.
We see that 
$$
\M_{l,l+1}(i,j) =
\begin{cases}
a_1 + \cdots + a_N & \text{ if } i= j=l+2,\\
b & \text{ if } 1 \le i= j \le l+1, \\
b & \text{ if } i + j = 2l +5, \, 1 \le i \le l+1, \\
c & \text{ if } l+3 \le i= j \le 2l+2, \\
c & \text{ if } i = 2l + 2, j= 2l +3, 2l+4, \\
0 & \text{ otherwise. } 
\end{cases}
$$  
Different from the symbolic matrix system for $\sccnrt$ is only the $l+2$-th row in 
$\M_{l,l+1}$.
The matrix $I_{l,l+1}$ is the same as the one for $\sccnrt$.
Let 
 $(M_{l,l+1},I_{l,l+1})_{l\in \Zp}$
be its nonnegative matrix system.
Hence we have 
$$
M_{l,l+1}^t(i,j) - I_{l,l+1}^t(i,j)  =
\begin{cases}
N & \text{ if } i=j=l+2,\\
1 & \text{ if } 2 \le i= j \le 2l+2, \, i\ne l+2,\\
1 & \text{ if } i + j = 2l +5, \, 1 \le j \le l+1, \\
-1 & \text{ if } 2 \le i = j+1\le 2l + 2, \\
0 & \text{ otherwise. } 
\end{cases}
$$  
By considering the kernels and cokernels of the following matrices
$B_{l,l+1}, l\in {\Bbb N}$ 
defined by 
$$
B_{l,l+1}(i,j)  =
\begin{cases}
N & \text{ if } i= j=l+2, \\
1 & \text{ if } 2 \le i= j \le 2l+2, \, i\ne l+2,\\
1 & \text{ if } (i, j) = (2l +4, 1), (2l+3, 2), \\
-1 & \text{ if } i=2, \, j=1, \\
0 & \text{ otherwise, } 
\end{cases}
$$  
that is 
$$
\setcounter{MaxMatrixCols}{16}
B_{l,l+1}
=
\begin{bmatrix}
            0&                         \hdotsfor{5}     & 0   & 0   & 0    & \hdotsfor{5}             \\
           -1& 1& 0    &               \hdotsfor{3}     & 0   & 0   & 0    & \hdotsfor{5}             \\
            0& 0& 1    & 0     &       \hdotsfor{2}     & 0   & 0   & 0    & \hdotsfor{5}             \\
\hdotsfor{1} & 0& 0    & 1     & 0    &\hdotsfor{1}     & 0   & 0   & 0    & \hdotsfor{5}             \\
\hdotsfor{2}    & \cdot& \cdot &\cdot &\cdot            &\cdot&\cdot& \cdot& \hdotsfor{5}             \\
\hdotsfor{3}           & \cdot &\cdot &\cdot            &\cdot&\cdot& \cdot& \hdotsfor{5}             \\
\hdotsfor{4}                   & 0    & 0               & 1   & 0   & 0    & \hdotsfor{5}             \\
\hdotsfor{5}                          & 0               & 0   & N   & 0    & \hdotsfor{5}             \\
\hdotsfor{6}                                            & 0   & 0   & 1    & 0   &              \hdotsfor{4} \\
\hdotsfor{5}                                      & 0   & 0   & 0   & 0    & 1   & 0    &       \hdotsfor{3} \\ 
\hdotsfor{4}                               & 0    & 0   & 0   & 0   & 0    & 0   & 1    & 0    &\hdotsfor{2} \\
\hdotsfor{3}           & \cdot & \cdot&\cdot      & \hdotsfor{3}    &\cdot &\cdot&\cdot &\cdot &\cdots       \\ 
\hdotsfor{2}    & \cdot& \cdot & \cdot&\hdotsfor{5}                        &\cdot&\cdot &\cdot &\cdot        \\
\hdotsfor{1} & 0& 0    &0      & \hdotsfor{7}                                    &0     &  0   &1            \\
            0& 1& 0    &         \hdotsfor{10}                                                 &0            \\
            1& 0&                \hdotsfor{11}                                                 &0 
 \end{bmatrix}. 
$$
We can similarly show that
$$
K_1(\sccnr )\cong 0 
$$
and
\begin{align*}
K_0( \sccnr ) 
= & \varinjlim 
\{ {\Bbb Z}^{2l + 4} / B_{l,l+1}{\Bbb Z}^{2l+2}, \overline{I^t}_{l,l+1} \} \\
= & \varinjlim \{ {\Bbb Z}/N{\Bbb Z} \oplus {\Bbb Z}\oplus{\Bbb Z}, L \} \\
\cong & {\Bbb Z}/N{\Bbb Z} \oplus {\Bbb Z}.
\end{align*}
Therefore we have
\begin{theorem}
\begin{align*}
K_0(\sccnr )& \cong  K_0(\sccnrt ) \cong {\Bbb Z}/N{\Bbb Z} \oplus {\Bbb Z},\\
K_1(\sccnr )& \cong  K_1(\sccnrt ) \cong 0.
\end{align*}
\end{theorem}

\begin{cor}
For $N,N'\in {\Bbb N},\ N \ne N'$,
$\sccnr$ and $sc({\mathcal C}^{(N')}_{reset})$
are not flow equivalent to each other. 
\end{cor}
\begin{proof}
K-groups are invariants of flow equivalence (\cite{Ma2001a}).
\end{proof}

\bibliographystyle{amsplain}

\end{document}